\documentclass[12pt]{article}
\usepackage[T1]{fontenc}
\usepackage[cp1250]{inputenc}
\usepackage[dvips]{graphicx}
\usepackage{amssymb}
\usepackage{amsmath}
\usepackage{latexsym}
\usepackage{mathabx}
\newcommand{\1}{{\rm 1\!\!1}}
\newcounter{myownsection}
\def\myownsection{\refstepcounter{myownsection} \setcounter{equation}{0}}

\begin{document}
\begin{center}
{\bf DECOMPOSITIONS OF THE FREE PRODUCT OF GRAPHS}
\footnote{This work is partially supported by the MNiSW research grant 1 P03A 013 30 and by 
the EU Network QP-Applications, Contract No. HPRN-CT-2002-00279}\\[60pt]
{\bf Luigi Accardi}\\[5pt]
{\it Centro Vito Volterra\\
Universita di Roma Tor Vergata\\
e-mail accardi@volterra.uniroma2.it}\\[10pt]
{\bf Romuald Lenczewski}\\[5pt]
{\it Institute of Mathematics and Computer Science\\
Wroc{\l}aw University of Technology\\
e-mail Romuald.Lenczewski@pwr.wroc.pl}\\[10pt]
{\bf Rafa{\l} Sa{\l}apata}\\[5pt]
{\it Institute of Mathematics and Computer Science\\
Wroc{\l}aw University of Technology\\
e-mail Rafal.Salapata@pwr.wroc.pl}\\[40pt]
\end{center}
\begin{abstract}
We study the free product of rooted graphs and its various decompositions 
using quantum probabilistic methods.
We show that the free product of rooted graphs is canonically associated with free independence, 
which completes the proof of the conjecture that there exists a product of rooted graphs canonically 
associated with each notion of noncommutative independence which arises in the axiomatic theory.
Using the {\it orthogonal product} of rooted graphs, we decompose the 
{\it branches} of the free product of rooted graphs 
as `alternating orthogonal products'.  This leads to alternating decompositions
of the free product itself, with the star product or the comb product followed by 
orthogonal products. These decompositions correspond to the recently studied
decompositions of the free additive convolution of probability measures in terms boolean and 
orthogonal convolutions, or monotone and orthogonal convolutions.
We also introduce a new type of {\it quantum decomposition} of the 
free product of graphs, where the distance partition of the set of vertices
is taken with respect to a set of vertices instead of a single vertex. 
We show that even in the case of widely studied graphs this yields new and 
more complete information on their spectral properties, like spectral
measures of a (usually infinite) set of cyclic vectors under the action 
of the adjacency matrix.
\\[10pt]
\end{abstract}
\newpage
\myownsection
\begin{center}
{\sc 1. Introduction}
\end{center}
Graph theory plays an important role in many branches of mathematics and its applications.
In particular, in solid state physics the idea that the study of a free dynamics on a non-homogeneous graph
might be equivalent to the study of an interacting dynamics on a homogeneous graph (a discrete
version of the basic idea of general relativity) found interesting applications in the study of
Bose--Einstein condensation [7,8,23].

In many of these applications, the main focus was not so much on the combinatorial properties of single graphs,
as on the analytical aspects of the asymptotics of large graphs (when the number of vertices tends to infinity).
In that connection, typical objects of interest are spectral distributions of the adjacency matrix with
respect to specially chosen states, or the full spectrum.

From these investigations an interesting class of graphs has emerged, namely those which are
built from subgraphs expressible as products of simpler graphs. Intuitively, a {\it product\/}
of graphs is a rule to construct a new graph by {\it glueing together\/} two given graphs subject to additional
conditions (like associativity). It is known that there exist several different products among graphs, including
those studied extensively in discrete mathematics like the lexicographic product, the Cartesian product or
the strong product.

On the other hand, the experience of quantum probability teaches us that, at an algebraic level,
certain types of products of quantum probability spaces correspond to different notions of 
stochastic independence [31,28]. Recall that in the concept of a quantum probability space one has to 
distinguish a state, which, in the simplest graph setting, corresponds to a distinguished vertex called {\it root}.
It is therefore natural to conjecture that certain types of products of rooted graphs
could be canonically associated with the main notions of stochastic independence. 
Moroever, one might expect that such graph products are important types of products, 
from which one could not only construct more complicated graphs, but also obtain information 
about their spectra, using (well-established or entirely new) quantum probabilistic techniques.

The well-known case of the Cayley graph of a free product of groups and its relation to
free independence of Voiculescu [33] and to the free product of states [3,33] 
can be viewed as the first evidence that such a conjecture might be true.
That such a relation holds also in the general case of the free product of rooted graphs
can be shown using the free probability techniques.
Thus, we explicitly state and prove the fact that the free product of rooted graphs, introduced by Znojko [39] for symmetric graphs and generalized by Quenell [30] and Gutkin [11] to rooted graphs, is canonically associated with 
the notion of free independence. More generally, we show that the hierarchy of {\it $m$--freeness\/}
introduced by one of the authors [18], corresponds to a natural hierarchy of
{\it $m$--free products\/} of graphs as well as to the natural inductive definition of the free product 
[39,11].

The next important evidence supporting this conjecture was the discovery, by Accardi, Ben Ghorbal
and Obata [1], that the {\it comb product\/} of rooted graphs 
is canonically related to the {\it monotone independence\/} [22,26].
A similar connection between the {\it star product\/} of
rooted graphs and {\it boolean independence} [32] was established by one of the authors [17] and Obata [29].
Since it is well-known that the Cartesian product of graphs is naturally related to
tensor (or boson) independence, the correspondence between certain types of products
of rooted graphs and the main notions of stochastic independence is completed.

In all the above-mentioned cases, the canonical relation between a notion of product of two graphs,
${\cal G}_{1}=(V_{1},E_{1})$ and ${\cal G}_{2}=(V_{2},E_{2})$,
and a notion of independence ${\cal I}$ is realized by showing that the adjacency matrix of their
product is naturally split into a sum of operator random variables which are
${\cal I}$--independent with respect to a given state.
In all cases this decomposition can be obtained by embedding the algebra of operators on the
$l^2$--space of the product graph into an appropriate tensor product, paralleling the
construction of [18,19].

Let us remark in this context that the Cartesian, comb and star products
appear to be `basic' graph products since the corresponding 
vertex sets are subsets of $V_{1}\times V_{2}$, whereas
the vertex set of the free product of graphs is the free product 
$V_{1}*V_{2}$ of rooted sets and the construction of ${\cal G}_{1}*{\cal G}_{2}$ involves
infinitely many copies of ${\cal G}_{1}$ and ${\cal G}_{2}$. 
This is the main reason why the free product of graphs has this peculiar feature 
that it admits a variety of natural decompositions. In particular, the decomposition
into the sum of `freely independent subgraphs', although the most natural 
from the point of view of free independence, is not always the most
intuitive or the most convenient. 

In particular, we find new decompositions related to the `growth' of ${\cal G}_{1}*{\cal G}_{2}$
exhibited by its inductive definitions. Motivated by the recent work of one of the authors on 
the decompositions of the free additive convolution of probability measures [21], 
we study a new type of `basic' 
product of rooted graphs called the {\it orthogonal product} of graphs, related to the 
{\it orthogonal convolution} of probability measures introduced in [21]. 
We show that this is the orthogonal product which is the main building block 
of ${\cal G}_{1}*{\cal G}_{2}$ since it is responsible for its `growth'.
In fact, it allows us to decompose its `branches' [30] into products of 
alternating  ${\cal G}_{1}$ and ${\cal G}_{2}$.
Since one obtains the free product by taking the comb product 
or the star product of `branches', as shown in [30] (we use quantum probabilistic techniques
to simplify the proofs presented there), we arrive at two alternating decompositions of the free product 
of graphs --
the {\it comb-orthogonal decomposition} and the {\it star-orthogonal decomposition}.
More importantly, using the orthogonal convolution, one can study
spectral distributions of free products of uniformly locally finite graphs in a very intuitive manner
and see their direct relation to continued fractions, especially mixed periodic
Jacobi continued fractions [14] (without using the R-transforms).

Finally, in order to get a more detailed information about the structure of the 
spectrum of ${\cal G}_{1}*{\cal G}_{2}$, we introduce a new type of {\it quantum decomposition}
of the adjacency matrix of a given graph ${\cal G}$. 
This decomposition is based on a new type of distance partition
$V=\bigcup_{n=0}^{\infty}{\cal V}_{n}$ of the set of vertices,
where ${\cal V}_{n}$ is the set of vertices of ${\cal G}$, whose distance from a
a set of vertices (instead of a single vertex) is equal to $n$. This leads to a different 
quantum decomposition of the adjacency matrix $A({\cal G})$ into the sum of a creation,
annihilation and diagonal operators. It allows us to derive a cyclic direct sum decomposition of
the Hilbert space $l_{2}(V)$ together with the spectral distributions associated with different 
cyclic (vacuum) vectors.

For classical (random walk) methods applied to the free products of Cayley graphs and other infinite 
graphs, we refer the reader to [2,6,9,10,15,16,37,38] and references 
contained there. \\[10pt]

\myownsection
\begin{center}
{\sc 2. Notation}
\end{center}
By a {\it rooted set} we understand a pair $(X,e)$, where $X$ is a countable
set and $e$ is a distinguished element of $X$ called {\it root}.
By a {\it rooted graph} we understand a pair $({\cal G},e)$, where
${\cal G}=(V,E)$ is a non-oriented graph with the set of {\it vertices} $V=V({\cal G})$,
and the set of edges $E=E({\cal G})\subseteq \{\{x,x'\}: \; x,x'\in V, x\neq x'\}$
and $e\in V$ is a distinguished vertex called the {\it root}.
We will also denote by ${\cal G}$ the rooted graph $({\cal G},e)$ if no confusion arises, 
especially if the graph
is {\it symmetric}, i.e. for any $x\neq x'$ there exists an authomorphism $\tau$ of
${\cal G}$ for which $\tau(x)=x'$ (in other words, all vertices are equivalent).

For rooted graphs we will use the notation
\begin{equation}
V^{0}=V\setminus \{e\}.
\end{equation}
Two vertices $x,x'\in V$ are called {\it adjacent} if $\{x,x'\}\in E$,  i.e.
vertices $x,x'$ are connected with an edge. Then we write $x \sim x'$.
Simple graphs have no loops, i.e. $\{x,x\}\notin E$ for all $x\in V$.
The {\it degree} of $x\in V$ is defined by $\kappa(x)=|\{x'\in V:x'\sim x\}|$,
where $|I|$ stands for the cardinality of $I$.
A graph is called {\it locally finite} if $\kappa(x)<\infty$ for every $x\in V$.
It is called {\it uniformly locally finite} if ${\rm sup}\{\kappa(x): x\in V\}<\infty$.

For $x\in V$, let $\delta(x)$ be the indicator function of the one-element set $\{x\}$.
Then $\{\delta(x),\,_x\in V\}$ is an orthonormal basis of the Hilbert space $l^{2}(V)$ of
square integrable functions on the set $V$, with the usual inner product.

The {\it adjacency matrix} $A=A({\cal G})$ of ${\cal G}$ is a 0-1 matrix defined by
\begin{equation}
A_{x,x'}=\left\{
\begin{array}{ll}
1 & {\rm if} \;\; x\sim x'\\
0 & {\rm otherwise}
\end{array}
\right.
\end{equation}
We identify $A$ with the densely defined symmetric operator on $l^{2}(V)$ defined by
\begin{equation}
A\delta(x)=\sum_{x\sim x'}\delta(x')
\end{equation}
for $x\in V$. Notice that the sum on the right-hand-side is finite since our graph is assumed
to be locally finite. It is known that $A({\cal G})$ is bounded iff ${\cal G}$ is uniformly
locally finite. If $A({\cal G})$ is essentially self-adjoint, its closure is called the
{\it adjacency operator} of ${\cal G}$ and its spectrum - the spectrum of ${\cal G}$.

The unital algebra generated by $A$, i.e. the algebra of polynomials in $A$, is called
the {\it adjacency algebra} of ${\cal G}$ and is denoted by ${\cal A}({\cal G})$ or simply
${\cal A}$.

In this paper by a graph we shall always
understand a non-oriented connected locally finite simple graph with
a non-empty set of edges.
Any rooted graph of type $({\cal G},e)$, where ${\cal G}$ is a graph
in this sense, will also be called a graph if no confusion arises.
\\[10pt]
\newpage
\myownsection
\begin{center}
{\sc 3. Convolutions, transforms and graph products}
\end{center}
By the {\it spectral distribution of $A({\cal G})$ in a state} $\psi$ on $l_{2}(V)$
we understand the measure $\mu$ for which
\begin{equation}
\psi (A^{n})= \int_{{\mathbb R}}x^{n}\mu(dx), \;\; n\in {\mathbb N} \cup \{0\}
\end{equation}
and by the {\it spectral distribution of the rooted graph $({\cal G},e)$}
we understand the spectral distribution of $A({\cal G})$ in the state
$\varphi _{e}(.)=\langle . \delta(e),\delta(e)\rangle$.
The spectral distribution of $({\cal G},e)$ is
important in the evaluation of the spectrum ${\rm spec}({\cal G})$ of the graph ${\cal G}$.
In some cases (homogenous trees and $n$-ary trees are the easiest examples)
it is even so that ${\rm spec}({\cal G})$
agrees with the support of the spectral distribution of $({\cal G},e)$.

For any probability distribution $\mu$ with the sequence of
moments $(M_{n})_{n\geq 0}$, we define its moment generating function as the formal
power series 
$$
M_{\mu}(z)=\sum_{n=0}^{\infty}M_{n}z^{n}.
$$ 
The corresponding Cauchy
transform, K-transform, reciprocal Cauchy transform and R-transform are defined respectively by the formal
power series
\begin{eqnarray}
G_{\mu}(z)&=&\frac{1}{z}M_{\mu}\left(\frac{1}{z}\right)\\
K_{\mu}(z)&=&z-\frac{1}{G_{\mu}(z)}\\
F_{\mu}(z)&=&\frac{1}{G_{\mu}(z)}\\
R_{\mu}(z)&=& -\frac{1}{z} +G_{\mu}^{-1}(z).
\end{eqnarray}

Let $({\cal G}_{1},e_{1})$ and $({\cal G}_{2},e_{2})$ be two graphs with adjacency matrices
$A_{1}=A({\cal G}_{1})$, $A_{2}=A({\cal G}_{2})$ and spectral distributions
$\mu$, $\nu$, respectively. By $\mu\uplus \nu$ we denote the
{\it boolean convolution} of $\mu$ and $\nu$ associated with boolean independence
[32]. By $\mu\vartriangleright \nu$ we denote
the {\it monotone convolution} associated with monotone independence [27].
By $\mu \boxplus \nu$ we denote the {\it free additive
convolution} associated with free independence [34].
Finally, by $\mu\vdash\nu$ we denote the {\it orthogonal convolution}, recently
introduced in [21].
The following identities hold:
\begin{eqnarray}
K_{\mu \uplus \nu}(z)&=&K_{\mu}(z)+K_{\nu}(z)\\
R_{\mu \,\boxplus\, \nu}(z)&=&R_{\mu}(z)+R_{\nu}(z)\\
F_{\mu \vartriangleright \nu}(z)&=&F_{\mu}(F_{\nu}(z)) \\
K_{\mu\,\vdash\,\nu}(z)&=& K_{\mu}(F_{\nu}(z)).
\end{eqnarray}
The above relations can be treated as definitions of the associated convolutions, 
although these are usually introduced by using some notion of noncommutative 
`independence' which parallels the connection between the usual (classical) convolution
of distributions (measures) and the notion of classical independence.
Note that the K- and R-transforms are additive under the considered convolutions
(see [32], [34], [37]), thus they play the role of the logarithm of the Fourier transform. 
In the case of the monotone and orthogonal convolutions, addition of transforms is replaced
by composition (see [27] and [21]).\\
\indent{\par}
{\sc Proposition 3.1.}
{\it The following relations hold:}
\begin{eqnarray}
F_{\mu\uplus \nu}(z) &=&F_{\mu}(z)+F_{\nu}(z)-z\\
F_{\mu \vdash \nu}(z) &=& F_{\mu}(F_{\nu}(z))-F_{\nu}(z)+z
\end{eqnarray}
\indent{\par}
{\it Proof.}
These are straightforward consequences of (3.3),(3.4) and (3.6),(3.9).\hfill $\blacksquare$\\
\indent{\par}
It is natural that the free additive convolution and the R-transform appear in the context
of spectral theory of free product graphs. It seems less so in the case of the other
three types of convolutions. However, we will show that one can decompose the free product
of graphs using the products of graphs associated with these convolutions.
Let us give the definitions of these products.   \\
\indent{\par}
{\sc Definition 3.1.}
The {\it comb product} of rooted graphs $({\cal G}_{1},e_1)$ and $({\cal G}_{2},e_2)$
is the rooted graph $({\cal G}_{1}\vartriangleright {\cal G}_{2},e)$
obtained by attaching a copy of ${\cal G}_{2}$ by its root $e_{2}$
to each vertex of ${\cal G}_{1}$, where we denote by $e$ the vertex obtained
by identifying $e_1$ and $e_2$. If no confusion arises, we denote the comb product by
${\cal G}_{1}\vartriangleright {\cal G}_{2}$. If we identify its set of vertices 
with $V_{1}\times V_{2}$, then its root is identified with $e_{1}\times e_{2}$. \\
\indent{\par}
Note that the comb product of rooted graphs is not commutative and it depends on the choice
of the root. Let us also remark that the definition given in [1] is equivalent to the one given
above, except that in our definition the information about the role of the root $e_2$ in the glueing
is encoded in the definition of the {\it rooted} graph $({\cal G}_{2},e_2)$.
Moreover, as our product is taken in the (natural) category of rooted graphs,
we define the root of the comb product to be $e$, which makes the product associative.\\
\indent{\par}
{\sc Theorem 3.2.} [1] 
{\it Let $({\cal G}_{1},e_{1})$ and $({\cal G}_{2},e_{2})$ be rooted graphs with
spectral distributions $\mu$ and $\nu$, respectively. Then, the adjacency matrix of 
their comb product can be decomposed as}
\begin{equation}
A({\cal G}_{1}\vartriangleright {\cal G}_{2}) =A^{(1)}+A^{(2)}
\end{equation}
{\it where $A^{(1)}$ and $A^{(2)}$ are monotone independent w.r.t. 
$\varphi(.)=\langle . \delta(e), \delta(e)\rangle $. Moreover, the spectral distribution of 
$({\cal G}_{1}\vartriangleright {\cal G}_{2},e)$ 
is given by $\mu \vartriangleright \nu$.}\\
\indent{\par}
{\sc Definition 3.2.}
The {\it star product} of $({\cal G}_{1},e_1)$ and $({\cal G}_{2},e_2)$ is the graph
$({\cal G}_{1}\star{\cal G}_{2},e)$ obtained by attaching a copy of
${\cal G}_{2}$ by its root $e_2$ to the root $e_1$ of ${\cal G}_{1}$,
where we denote by $e$ the vertex obtained by
identifying $e_1$ and $e_2$. If no confusion arises,
we also denote the star product by ${\cal G}_{1}\star{\cal G}_{2}$.
If we identify its set of vertices with  $V_{1}\star V_{2}:=
(V_{1}\times \{e_2\}) \cup (\{e_1\}\times V_{2})$, then its root is identified
with $e_{1}\times e_{2}$.\\
\indent{\par}
{\sc Theorem 3.3.} [17,29] 
{\it Let $({\cal G}_{1},e_{1})$ and $({\cal G}_{2},e_{2})$ be rooted graphs with
spectral distributions $\mu$ and $\nu$, respectively. Then, the adjacency matrix of their
star product can be decomposed as}
\begin{equation}
A({\cal G}_{1}\star {\cal G}_{2}) =A^{(1)}+A^{(2)}
\end{equation}
{\it where $A^{(1)}$ and $A^{(2)}$ are boolean independent w.r.t. $\varphi$,
where $\varphi(.)=\langle .\delta(e), \delta(e)\rangle$.
Moreover, the spectral distribution of $({\cal G}_{1}\star {\cal G}_{2},e)$ 
is given by $\mu \uplus \nu$.}\\[10pt]

\myownsection
\begin{center}
{\sc 4. Orthogonal product of graphs}
\end{center}
Let us introduce now a new basic product of rooted graphs called `orthogonal', which is
related to the orthogonal convolution introduced in [21]. Using this new product, together with the 
comb product of graphs (or, the star product of graphs), one can construct their free 
product ${\cal G}_{1}*{\cal G}_{2}$ using copies of ${\cal G}_{1}$ and ${\cal G}_{2}$.
This is an application of the more general theory of constructing the free additive 
convolution from the orthogonal and monotone (or, orthogonal and boolean) 
convolutions given in [21] and we use the results contained there. \\
\indent{\par}
{\sc Definition 4.1.}
The {\it orthogonal product} of two rooted graphs $({\cal G}_{1},e_{1})$ and $({\cal G}_{2},e_{2})$ 
is the rooted graph $({\cal G}_{1}\vdash {\cal G}_{2},e)$ obtained by attaching a copy of ${\cal G}_{2}$ by its root $e_{2}$ to each vertex of ${\cal G}_{1}$ but the root $e_{1}$, where $e$ is taken to be equal to $e_{1}$.
If its set of vertices is identified with $V_{1}\vdash V_{2}:=(V_{1}^{0}\times V_{2})\cup \{e_{1}\times e_{2}\}$
then $e$ is identified with $e_{1}\times e_{2}$.\\
\indent{\par}
It is worth noting that the orthogonal product of graphs resembles their comb product. 
The difference is that in the comb product the second graph is glued by its root to all 
vertices of the first graph, whereas in the orthogonal product the second graph is glued 
to all vertices {\it but the root} of the first graph. An example of the orthogonal product 
of graphs is given in Fig.1.\\
\indent{\par}
{\sc Example 4.1.}

\unitlength=1mm
\special{em.linewidth 2pt}
\linethickness{0.5pt}
\begin{picture}(60.00,30.00)(-30.00,7.00)
\put(0.00,10.00){\line(1,0){20.00}}
\put(35.00,10.00){\line(0,1){6.00}}
\put(35.00,16.00){\line(1,2){3.00}}
\put(35.00,16.00){\line(-1,2){3.00}}

\put(0.00,10.00){\circle*{1.00}}
\put(10.00,10.00){\circle*{1.00}}
\put(20.00,10.00){\circle*{1.00}}

\put(35.00,10.00){\circle*{1.00}}
\put(35.00,16.00){\circle*{1.00}}
\put(32.00,22.00){\circle*{1.00}}
\put(38.00,22.00){\circle*{1.00}}

\put(-2.00, 6.00){\footnotesize $e_{1}$}
\put(8.00,6.00){\footnotesize $x$}
\put(18.00,6.00){\footnotesize $x'$}
\put(36.00,8.00){\footnotesize $e_{2}$}
\put(36.00,14.00){\footnotesize $y$}
\put(40.00,22.00){\footnotesize $y''$}
\put(28.00,22.00){\footnotesize $y'$}
\put(8.00,14.00){\footnotesize ${\cal G}_{1}$}
\put(26.00,14.00){\footnotesize ${\cal G}_{2}$}
\put(55.00,14.00){\footnotesize ${\cal G}_{1}\vdash {\cal G}_{2}$}

\put(70.00,10.00){\line(1,0){20.00}}
\put(80.00,10.00){\line(0,1){6.00}}
\put(90.00,10.00){\line(0,1){6.00}}
\put(80.00,16.00){\line(-1,2){3.00}}
\put(80.00,16.00){\line(1,2){3.00}}
\put(90.00,16.00){\line(-1,2){3.00}}
\put(90.00,16.00){\line(1,2){3.00}}

\put(70.00,10.00){\circle*{1.00}}
\put(80.00,10.00){\circle*{1.00}}
\put(90.00,10.00){\circle*{1.00}}

\put(80.00,16.00){\circle*{1.00}}
\put(90.00,16.00){\circle*{1.00}}
\put(83.00,22.00){\circle*{1.00}}
\put(77.00,22.00){\circle*{1.00}}
\put(93.00,22.00){\circle*{1.00}}
\put(87.00,22.00){\circle*{1.00}}

\end{picture}
\\
\begin{center}
{\it Fig. 1.} Orthogonal product ${\cal G}_{1}\vdash {\cal G}_{2}$\\[5pt]
\end{center}
\indent{\par}
The notion of the orthogonal product of graphs is related to the concept
of {\it orthogonal subalgebras} introduced in [21]. Note that the concept of `orthogonality' 
involves two functionals (states), but it is quite different from `conditional freeness'[5].   \\
\indent{\par}
{\sc Definition 4.2.}
Let $({\cal A},\varphi, \psi)$ be a unital algebra with a pair
of linear normalized functionals and let ${\cal A}_{1}$ and ${\cal A}_{2}$
be non-unital subalgebras of ${\cal A}$.
We say that ${\cal A}_{2}$ is {\it orthogonal} to ${\cal A}_{1}$
with respect to $(\varphi, \psi)$ if\\
\indent{\par}
(i)
$\;\varphi(bw_{2})=\varphi(w_{1}b)=0$
\indent{\par}
(ii)
$\varphi(w_{1}a_{1}ba_{2}w_{2})=\psi(b)
\left(\varphi(w_{1}a_{1}a_{2}w_{2})- \varphi(w_{1}a_{1})\varphi(a_{2}w_{2})\right)$\\[5pt]
for any $a_{1},a_{2}\in {\cal A}_{1}$, $b\in {\cal A}_{2}$
and any elements $w,v$ of the unital algebra ${\rm alg}({\cal A}_{1},{\cal A}_{2})$
generated by ${\cal A}_{1}$ and ${\cal A}_{2}$.
We say that the pair $(a,b)$ of elements of ${\cal A}$ is {\it orthogonal} with respect to $(\varphi, \psi)$
if the algebra generated by $a\in {\cal A}$ is orthogonal to the algebra generated by
$b\in {\cal A}$ .\\
\indent{\par}
In analogy to Theorems 3.2-3.3, one can decompose the adjacency matrix of the 
orthogonal product of graphs. The proof is based on the tensor product realization of
orthogonal subalgebras [21]. Note that tensor product realizations of noncommutative 
random variables, originated in [18] for boolean, $m$-free and free random variables,
(see also the review paper [20]), are especially useful in the context of graph products
since the projections introduced in this scheme tell us how the graphs should be glued
together. This technique was later used in a number of papers ([1],[12], [19], [29]). 
\\
\indent{\par}
{\sc Theorem  4.1.}
{\it Let $({\cal G}_{1},e_{1})$ and $({\cal G}_{2},e_{2})$ be rooted graphs with
spectral distributions $\mu$ and $\nu$, respectively.
Then, the adjacency matrix of their orthogonal product can be decomposed as}
\begin{equation}
A({\cal G}_{1}\vdash {\cal G}_{2}) =A^{(1)}+A^{(2)}
\end{equation}
{\it where the pair $(A^{(1)},A^{(2)})$ is orthogonal w.r.t. $(\varphi, \psi)$, where
$\varphi$ and $\psi$ are states associated with 
vectors $\delta(e), \delta(v)\in l_{2}(V_{1}\vdash V_{2})$ and $v\in V_{1}^{0}$.
Moreover, the spectral distribution of $({\cal G}_{1}\vdash {\cal G}_{2},e)$ 
is given by $\mu \vdash \nu$.}\\
\indent{\par}
{\it Proof.}
In order to prove the decomposition, it is convenient
to identify the adjacency matrix of ${\cal G}_{1}\vdash {\cal G}_{2}$ with the sum
$$
A({\cal G}_{1}\vdash {\cal G}_{2})=A_{1}\otimes P_{\xi_{2}}+P_{\xi_{1}}^{\perp} \otimes A_{2}
$$
on the Hilbert space $l_{2}(V_{1}\vdash V_{2})\subset l_{2}(V_{1}\times V_{2})\cong l_{2}(V_{1})\otimes l_{2}(V_{2})$, where $A_{i}$ is the adjacency matrix of ${\cal G}_{i}$ and $P_{\xi_{1}}^{\perp}=1-P_{\xi_{1}}$,
with $P_{\xi_{i}}$ denoting the projection onto ${\mathbb C}\xi_{i}$, where $\xi_{i}=\delta(e_{i})$
and $i=1,2$. Projection $P_{\xi_{1}}^{\perp}$  indicates that graph ${\cal G}_{2}$ 
should be glued to all vertices of ${\cal G}_{1}$ but the root, whereas projection $P_{\xi_{2}}$
indicates that graph ${\cal G}_{1}$ should be glued only to vertex $e_{2}$ of ${\cal G}_{2}$,
which reproduces Definition 4.1. It remains to take $\varphi$ and $\psi$ to be the states associated
with unit vectors  $\delta(e_{1})\times \delta(e_{2})$ and $\delta(v)\times \delta(e_{2})$, 
respectively, where $v$ is an arbitrary vertex from $V_{1}^{0}$.
Now, in view of [Theorem 4.1, 21], the above summands form a pair of orthogonal 
elements (of the algebra they generate) w.r.t. the pair of states $(\varphi, \psi)$.
The proof consists in checking (i)-(ii) of Definition 4.2 and in the case of graphs
is very similar to the general case.

It follows from [Corollary 4.2, 21] that the spectral distribution of such a sum in the state $\varphi$ 
is equal to the orthogonal convolution $\mu \vdash \nu$.
Nevertheless, we choose to show this fact here since we can present a proof
which nicely exhibits the relation between the 
comb product and the orthogonal product.
Namely, recall that ${\cal G}_{1}\vdash {\cal G}_{2}$ differs from ${\cal G}_{1}\vartriangleright {\cal G}_{2}$
by the fact that no copy of ${\cal G}_{2}$ is glued to the vertex $e_{1}$. Therefore, one can obtain
${\cal G}_{1}\vartriangleright {\cal G}_{2}$ by glueing ${\cal G}_{1}\vdash {\cal G}_{2}$ and
${\cal G}_{2}$ at their roots, which corresponds to their star product. Therefore,
$$
{\cal G}_{1}\vartriangleright {\cal G}_{2}=({\cal G}_{1}\vdash {\cal G}_{2})\star {\cal G}_{2}
$$
which leads to the following formula for their spectral distributions
$$
\mu\vartriangleright \nu =\sigma \uplus \nu
$$
where $\sigma $ is the spectral distribution of ${\cal G}_{1}\vdash {\cal G}_{2}$. 
Using transforms, we get
$$
F_{\mu}(F_{\nu}(z))=F_{\sigma}(z)+F_{\nu}(z)-z
$$
which, in view of the second equation of Proposition 3.1, gives our assertion.\hfill $\blacksquare$\\
\indent{\par}
{\sc Example 4.1.}
Let us apply Theorem 4.1 to the orthogonal product in Fig.1.
We have 
$$
G_{\mu}(z)=\frac{z^{2}-1}{z(z^{2}-2)}, \;\;\;
G_{\nu}(z)=\frac{z^{2}-2}{z(z^{2}-3)}
$$
and therefore, 
$$
K_{\mu}(z)=\cfrac{1}{z-\cfrac{1}{z}}, \;\;\;
F_{\nu}(z)=z-\cfrac{1}{z-\cfrac{2}{z}}.
$$
In view of Theorem 4.1 and (3.9), we obtain the explicit formula
for the Cauchy transform
$$
G_{\mu\, \vdash \,\nu}(z)=\frac{1}{z-K_{\mu}(F_{\nu}(z))}.
$$
Algebraic calculations lead to the continued fraction representation 
of $G_{\mu\,\vdash\,\nu}(z)$ associated with the sequences of Jacobi coefficients $\omega=(\omega_{n})=(1,2,3/2,5/6,4/15,12/5,0,\ldots)$ and $\alpha=(\alpha_{n})=(0,0, \ldots)$.
The corresponding measure is a discrete measure consisting of 7 atoms (since their 
explicit values and corresponding masses are rather complicated, we do not give them here).\\[10pt]

\myownsection
\begin{center}
{\sc 5. Free and $m$-free products of graphs}
\end{center}
In this Section we recall after [39] and [11] the definition
of the free product of rooted graphs and define a
corresponding approximating sequence of $m$-free products.

Consider rooted graphs $({\cal G}_{i},e_i)=(V_{i},E_{i},e_i)$,
where $i \in I$ and $I$ is a finite index set, and denote $V_{i}^{0}=V_i\setminus \{e_i\}$.
By the {\it free product of rooted sets} $(V_{i},e_i)$, $i\in I$, we understand
the rooted set  $(*_{i\in I}V_{i}, e)$, where
$$
*_{i\in I}V_{i}=
\{e\} \cup \{v_{1}v_{2}\ldots v_{m};\;v_{k}\in V_{i_{k}}^{0}\;\; {\rm and}\;\;
i_{1}\neq i_{2}\neq \ldots \neq i_{n},\; m\in {\mathbb N}\}
$$
and $e $ is the empty word.
For notational convenience, we will sometimes use words containing roots $e_k$
but then we shall always understand that $we_k=e_kw\equiv w$ where
$w\in *_{i\in I}V_{i}$, thus any $e_{k}$ will be treated as the `unit' or the empty word.
We are ready to give the definition of the free product of graphs.
\\
\indent{\par}
{\sc Definition 5.1.}
By the {\it free product of rooted graphs}
$({\cal G}_{i},e_i)$, $i \in I$,
we understand the rooted graph $(*_{i\in I}{\cal G}_{i}, e)$
with the set of vertices $*_{i\in I}V_{i}$ and the set of edges $*_{i\in I}E_{i}$
consisting of pairs of vertices from $*_{i\in I}V_{i}$ of the form
$$
*_{i\in I}E_{i}=
\{\{vu,v'u\}:\; \{v,v'\}\in \bigcup_{i\in I}E_{i}\;{\rm and}\;
u,vu,v'u\in *_{i\in I}V_{i}\}.
$$
We denote this product by
$*_{i\in I}({\cal G}_{i},e_i)$ or simply
$*_{i\in I}{\cal G}_{i}$ if no confusion arises.\\
\indent{\par}
Observe that if one of the vertices, $v$ or $v'$, is a root $e_{j}$, then
$\{vu, v'u\}\in \bigcup_{i\in I}E_{i}$ provided the other one forms an edge
with it -- we use here the convention that every $e_j$ can be treated as the 'unit'.
Let us also mention that the free product of graphs is commutative and associative (cf. [39]), 
which also follows from commutativity and associativity of the free product of states in view of Theorem 6.2.
Moreover, it is clear that the free product of a finite number of (uniformly) 
locally finite graphs is (uniformly) locally finite.

The most intuitive construction of the free product of graphs is given by some inductive
procedure which gives a sequence of growing graphs whose inductive limit is the free product of graphs. 
In fact, one natural procedure was given in [39], where
it was one of the equivalent definitions of the free product of graphs.
Interestingly enough, this procedure gives a sequence of iterates indexed by
$m\in {\mathbb N}$ which corresponds to the $m$-free product of states introduced in [18]. 
This leads us to the following formal definition
(an example is given in Fig.2).\\
\indent{\par}
{\sc Definition 5.2.}
By the $m$-{\it free product} of rooted graphs
$({\cal G}_{i},e_i)$, $i\in I$, we understand the subgraph
$(*_{i\in I}^{(m)}{\cal G}_{i}, e)$ of $(*_{i\in I}{\cal G}_{i}, e)$
obtained by restricting the set of vertices to words $w$ of lenght $|w|\leq m$.\\
\indent{\par}
{\sc Proposition 5.1.}
{\it It holds that}
$$
*_{i\in I}{\cal G}_{i}=\bigcup_{m=1}^{\infty}*_{i\in I}^{m}{\cal G}_{i}
$$
\indent{\par}
{\it Proof.}
This is an immediate consequence of Definitions 5.1-5.2.\hfill $\blacksquare$\\
\indent{\par}
{\sc Example 5.1.}
Consider two `segments' ${\cal G}_{1}\cong {\mathbb Z}_{2}$ and
${\cal G}_{2}\cong {\mathbb Z}_{2}$. They are graphs consisting of one edge
$x\sim e_{1}$ and $y\sim e_2$. Then
$$
V_{1}*V_{2}=\{e, x,y,xy,yx,xyx,yxy, \ldots \}
$$
$$
E_{1}*E_{2}=\{\{e , x\}, \{e ,y\}, \{x,yx\},\{y,xy\},\{yx,xyx\}, \{xy, yxy\}, \ldots \}
$$
and it is easy to see that
${\cal G}_{1}*{\cal G}_{2}\cong {\mathbb Z}$, where ${\mathbb Z}$
denotes the two-way infinite path (or, 1-dimensional integer lattice) with the root at $0$.\\
\indent{\par}
{\sc Example 5.2.}
Let ${\cal G}_{1}=(V_{1},E_{1})$ and ${\cal G}_{2}=(V_{2},E_{2})$, where
$V_{1}=\{e_1,x,x'\}$, $V_{2}=\{e_2,y,y'\}$ and $e_1\sim x$, $x\sim x'$,
$e_{2}\sim y, y\sim y'$. Then we obtain
\begin{eqnarray*}
E_{1}*E_{2}&=&\{\{e,x\}, \{e ,y\}, \{x,x'\}, \{y,y'\}, \{x,yx\}, \{x',yx'\}\\
&& \{y,xy\}, \{y',xy'\},\{yx,y'x\},\{yx',y'x'\}, \{xy,x'y\},\\
&& \{xy',x'y'\}, \{yx,xyx\}, \{y'x,xy'x\}, \{xy, yxy\}, \{xy',yxy'\}, \ldots \}.
\end{eqnarray*}
In Fig.2 we draw the subgraph of this free product given by the
$4$-free product of $({\cal G}_{1},e_1)$ and $({\cal G}_{2},e_2)$.
Copies of ${\cal G}_{1}$ (drawn horizontally)
and ${\cal G}_{2}$ (drawn vertically) form a covering [39] of the product graph.
Using the rules of orthogonal glueing,
we label selected vertices by suitable words to show how to form the product graph.\\

\unitlength=1mm
\special{em.linewidth 2pt}
\linethickness{0.5pt}
\begin{picture}(100.00,65.00)(-35.00,-10.00)
\put(0.00,50.00){\line(1,0){50.00}}
\put(0.00,50.00){\line(0,-1){50.00}}

\put(-2.00, 51.00){\footnotesize $e$}
\put(23.00,51.00){\footnotesize $x$}
\put(48.00,51.00){\footnotesize $x'$}
\put(20.00,41.00){\footnotesize $yx$}
\put(20.00,32.00){\footnotesize $y'x$}
\put(30.00,44.00){\footnotesize $xyx$}
\put(36.00,42.00){\footnotesize $x'yx$}
\put(-4.00,24.00){\footnotesize $y$}
\put(-4.00,1.00){\footnotesize $y'$}
\put(7.00,1.00){\footnotesize $xy'$}
\put(17.00,1.00){\footnotesize $x'y'$}
\put(00.00,50.00){\circle*{1.00}}
\put(25.00,50.00){\circle*{1.00}}
\put(50.00,50.00){\circle*{1.00}}
\put(0.00,25.00){\circle*{1.00}}
\put(0.00,0.00){\circle*{1.00}}


\put(0.00,25.00){\line(1,0){15.00}}

\put(7.50,25.00){\circle*{1.00}}
\put(15.00,25.00){\circle*{1.00}}

\put(0.00,0.00){\line(1,0){15.00}}
\put(7.50,0.00){\circle*{1.00}}
\put(15.00,0.00){\circle*{1.00}}
\put(7.50,0.00){\line(0,-1){10.00}}
\put(15.00,0.00){\line(0,-1){10.00}}
\put(15.00,-10.00){\circle*{1.00}}
\put(7.50,-5.00){\circle*{1.00}}
\put(7.50,-10.00){\circle*{1.00}}
\put(15.00,-5.00){\circle*{1.00}}
\put(7.50,-5.00){\line(1,0){5.00}}
\put(15.00,-5.00){\line(1,0){5.00}}
\put(7.50,-10.00){\line(1,0){5.00}}
\put(15.00,-10.00){\line(1,0){5.00}}
\put(10.00,-5.00){\circle*{1.00}}
\put(12.50,-5.00){\circle*{1.00}}
\put(17.50,-5.00){\circle*{1.00}}
\put(20.00,-5.00){\circle*{1.00}}
\put(10.00,-10.00){\circle*{1.00}}
\put(12.50,-10.00){\circle*{1.00}}
\put(17.50,-10.00){\circle*{1.00}}
\put(20.00,-10.00){\circle*{1.00}}

\put(0.00,25.00){\line(1,0){15.00}}
\put(7.50,25.00){\circle*{1.00}}
\put(15.00,25.00){\circle*{1.00}}
\put(7.50,25.00){\line(0,-1){10.00}}
\put(15.00,25.00){\line(0,-1){10.00}}
\put(15.00,15.00){\circle*{1.00}}
\put(7.50,20.00){\circle*{1.00}}
\put(7.50,15.00){\circle*{1.00}}
\put(15.00,20.00){\circle*{1.00}}
\put(7.50,20.00){\line(1,0){5.00}}
\put(15.00,20.00){\line(1,0){5.00}}
\put(7.50,15.00){\line(1,0){5.00}}
\put(15.00,15.00){\line(1,0){5.00}}
\put(10.00,20.00){\circle*{1.00}}
\put(12.50,20.00){\circle*{1.00}}
\put(17.50,20.00){\circle*{1.00}}
\put(20.00,20.00){\circle*{1.00}}
\put(10.00,15.00){\circle*{1.00}}
\put(12.50,15.00){\circle*{1.00}}
\put(17.50,15.00){\circle*{1.00}}
\put(20.00,15.00){\circle*{1.00}}

\put(25.00,50.00){\line(0,-1){15.00}}
\put(50.00,50.00){\line(0,-1){15.00}}
\put(25.00,42.50){\line(1,0){10.00}}
\put(25.00,35.00){\line(1,0){10.00}}

\put(30.00,35.00){\line(0,-1){5.00}}
\put(30.00,42.50){\line(0,-1){5.00}}
\put(35.00,35.00){\line(0,-1){5.00}}
\put(35.00,42.50){\line(0,-1){5.00}}

\put(25.00,35.00){\circle*{1.00}}
\put(25.00,42.50){\circle*{1.00}}
\put(30.00,30.00){\circle*{1.00}}
\put(30.00,42.50){\circle*{1.00}}
\put(35.00,37.50){\circle*{1.00}}
\put(35.00,42.50){\circle*{1.00}}

\put(30.00,37.50){\circle*{1.00}}
\put(30.00,40.00){\circle*{1.00}}
\put(30.00,35.00){\circle*{1.00}}
\put(30.00,32.50){\circle*{1.00}}
\put(35.00,40.00){\circle*{1.00}}
\put(35.00,30.00){\circle*{1.00}}
\put(35.00,35.00){\circle*{1.00}}
\put(35.00,32.50){\circle*{1.00}}

\put(50.00,50.00){\line(0,-1){15.00}}
\put(50.00,42.50){\line(1,0){10.00}}
\put(50.00,35.00){\line(1,0){10.00}}

\put(55.00,35.00){\line(0,-1){5.00}}
\put(55.00,42.50){\line(0,-1){5.00}}
\put(60.00,35.00){\line(0,-1){5.00}}
\put(60.00,42.50){\line(0,-1){5.00}}

\put(50.00,35.00){\circle*{1.00}}
\put(50.00,42.50){\circle*{1.00}}
\put(55.00,37.50){\circle*{1.00}}
\put(55.00,42.50){\circle*{1.00}}
\put(60.00,37.50){\circle*{1.00}}
\put(60.00,42.50){\circle*{1.00}}

\put(55.00,40.00){\circle*{1.00}}
\put(55.00,30.00){\circle*{1.00}}
\put(55.00,35.00){\circle*{1.00}}
\put(55.00,32.50){\circle*{1.00}}
\put(60.00,40.00){\circle*{1.00}}
\put(60.00,30.00){\circle*{1.00}}
\put(60.00,32.50){\circle*{1.00}}
\put(60.00,35.00){\circle*{1.00}}

\put(33.00,0.00){${\cal G}_{1}$}
\put(40.00,0.00){\line(1,0){15.00}}
\put(40.00,0.00){\circle*{1.00}}
\put(38.00,-3.00){\footnotesize $e_1$}
\put(45.00,-3.00){\footnotesize $x$}
\put(47.50,0.00){\circle*{1.00}}
\put(55.00,0.00){\circle*{1.00}}
\put(52.00,-3.00){\footnotesize $x'$}
\put(60.00,0.00){${\cal G}_{2}$}
\put(65.00,-7.50){\line(0,1){15.00}}
\put(65.00,-7.50){\circle*{1.00}}
\put(65.00,0.00){\circle*{1.00}}
\put(67.00,-8.50){\footnotesize $y'$}
\put(67.00,-1.00){\footnotesize $y$}
\put(65.00,7.50){\circle*{1.00}}
\put(67.00,6.50){\footnotesize $e_2$}

\end{picture}
\\
\begin{center}
{\it Fig. 2.} 4-free product $({\cal G}_{1},e_1)*^{(4)}({\cal G}_{2},e_2)$
with selected vertices labelled.\\[10pt]
\end{center}

\indent{\par}
Notice that $({\cal G}_{1},e_1)*({\cal G}_{2},e_2)$ is, in general, not isomorphic
to $({\cal G}_{1},x)*({\cal G}_{2},y)$ if $e_1\neq x$ or $e_2 \neq y$
(in that case at least one graph should not be symmetric).\\
\indent{\par}
{\sc Example 5.3.}
Let ${\mathbb T}_n$ be the $n$-{\it ary rooted tree} with the root denoted
$e$, where $n\geq 2$.
First, observe that ${\mathbb T}_n$ is the
$n$-fold free power of ${\mathbb T}_1\cong {\mathbb Z}_{+}$,
\begin{equation}
{\mathbb T}_n\cong {\mathbb T}_1*{\mathbb T}_1* \ldots * {\mathbb T}_1 \;\; (n \;{\rm times})
\end{equation}
where ${\mathbb Z}_{+}=(V,E)$, with $V=\{0,1,2,\ldots\}$ and $E=\{\{x,x+1\},x\geq 0\}$.
In particular, the binary rooted tree is the free product of two copies of ${\mathbb T}_1$
(see Fig. 3).\\
\unitlength=1mm
\special{em.linewidth 0.5pt}
\linethickness{0.5pt}
\begin{picture}(140.00,60.00)(-10.00,5.00)

\put(10.00,10.00){\line(1,2){5.00}}
\put(20.00,10.00){\line(-1,2){5.00}}
\put(30.00,10.00){\line(1,2){5.00}}
\put(40.00,10.00){\line(-1,2){5.00}}

\put(50.00,10.00){\line(1,2){5.00}}
\put(60.00,10.00){\line(-1,2){5.00}}
\put(70.00,10.00){\line(1,2){5.00}}
\put(80.00,10.00){\line(-1,2){5.00}}

\put(15.00,20.00){\line(2,3){10.00}}
\put(35.00,20.00){\line(-2,3){10.00}}
\put(55.00,20.00){\line(2,3){10.00}}
\put(75.00,20.00){\line(-2,3){10.00}}

\put(25.00,35.00){\line(4,3){20.00}}
\put(65.00,35.00){\line(-4,3){20.00}}

\put(10.00,10.00){\circle*{1.50}}
\put(20.00,10.00){\circle*{1.50}}
\put(30.00,10.00){\circle*{1.50}}
\put(40.00,10.00){\circle*{1.50}}
\put(50.00,10.00){\circle*{1.50}}
\put(60.00,10.00){\circle*{1.50}}
\put(70.00,10.00){\circle*{1.50}}
\put(80.00,10.00){\circle*{1.50}}

\put(10.00,10.00){\line(-1,-2){2.00}}
\put(10.00,10.00){\line(1,-2){2.00}}
\put(20.00,10.00){\line(-1,-2){2.00}}
\put(20.00,10.00){\line(1,-2){2.00}}
\put(30.00,10.00){\line(-1,-2){2.00}}
\put(30.00,10.00){\line(1,-2){2.00}}
\put(40.00,10.00){\line(-1,-2){2.00}}
\put(40.00,10.00){\line(1,-2){2.00}}

\put(50.00,10.00){\line(-1,-2){2.00}}
\put(50.00,10.00){\line(1,-2){2.00}}
\put(60.00,10.00){\line(-1,-2){2.00}}
\put(60.00,10.00){\line(1,-2){2.00}}
\put(70.00,10.00){\line(-1,-2){2.00}}
\put(70.00,10.00){\line(1,-2){2.00}}
\put(80.00,10.00){\line(-1,-2){2.00}}
\put(80.00,10.00){\line(1,-2){2.00}}

\put(15.00,20.00){\circle*{1.50}}
\put(35.00,20.00){\circle*{1.50}}
\put(55.00,20.00){\circle*{1.50}}
\put(75.00,20.00){\circle*{1.50}}

\put(25.00,35.00){\circle*{1.50}}
\put(65.00,35.00){\circle*{1.50}}

\put(45.00,50.00){\circle*{1.50}}

\put(5.00,10.00){\footnotesize $x''$}
\put(14.00,10.00){\footnotesize $yx'$}
\put(23.50,10.00){\footnotesize $xyx$}
\put(33.00,10.00){\footnotesize $y'x$}
\put(51.00,10.00){\footnotesize $x'y$}
\put(61.00,10.00){\footnotesize $yxy$}
\put(72.00,10.00){\footnotesize $xy'$}
\put(81.00,10.00){\footnotesize $y''$}

\put(10.00,20.00){\footnotesize $x'$}
\put(37.00,20.00){\footnotesize $yx$}
\put(47.00,20.00){\footnotesize $xy$}
\put(77.00,20.00){\footnotesize $y'$}

\put(20.00,35.00){\footnotesize $x$}
\put(67.00,35.00){\footnotesize $y$}
\put(44.00,52.00){\footnotesize $e$}

\put(100.00,10.00){\line(1,0){18.00}}
\put(100.00,50.00){\line(0,-1){18.00}}

\put(100.00,10.00){\circle*{1.50}}
\put(105.00,10.00){\circle*{1.50}}
\put(110.00,10.00){\circle*{1.50}}
\put(115.00,10.00){\circle*{1.50}}

\put(100.00,50.00){\circle*{1.50}}
\put(100.00,45.00){\circle*{1.50}}
\put(100.00,40.00){\circle*{1.50}}
\put(100.00,35.00){\circle*{1.50}}

\put(99.00,6.00){\footnotesize $e_1$}
\put(104.00,6.00){\footnotesize $x$}
\put(109.00,6.00){\footnotesize $x'$}
\put(114.00,6.00){\footnotesize $x''$}
\put(96.00,49.00){\footnotesize $e_2$}
\put(96.00,44.00){\footnotesize $y$}
\put(96.00,39.00){\footnotesize $y'$}
\put(96.00,34.00){\footnotesize $y''$}

\put(100.00,15.00){\footnotesize ${\cal G}_{1}\cong {\mathbb T}_1$}
\put(105.00,40.00){\footnotesize ${\cal G}_{2}\cong {\mathbb T}_1$}

\end{picture}\\
\begin{center}
{\it Fig. 3.}
Binary tree  ${\mathbb T}_2\cong {\mathbb T}_1*{\mathbb T}_1$\\[5pt]
\end{center}
\indent{\par}
{\sc Example 5.4.}
Let ${\mathbb H}_n$ denote the {\it homogenous tree of order} $n$, where $n\geq 1$, with
the root $e$.
If $n$ is even, then ${\mathbb H}_n$ is the Cayley graph of a free group.
It is easy to see that
\begin{equation}
{\mathbb H}_n\cong {\mathbb Z}_{2}*{\mathbb Z}_{2}*\ldots {\mathbb Z}_{2}\;\; (n \;{\rm times})
\end{equation}
with the 'segment' ${\mathbb Z}_{2}$ as in Example 5.1.\\

\unitlength=1mm
\special{em.linewidth 0.5pt}
\linethickness{0.5pt}
\begin{picture}(100.00,60.00)(-30.00,0.00)

\put(-2.00,25.00){\line(1,0){54.00}}
\put(0.00,22.50){\line(1,0){5.00}}
\put(0.00,27.50){\line(1,0){5.00}}
\put(4.00,17.50){\line(1,0){10.00}}
\put(4.00,32.50){\line(1,0){10.00}}
\put(13.00,10.00){\line(1,0){24.00}}

\put(15.00,7.50){\line(1,0){5.00}}
\put(15.00,12.50){\line(1,0){5.00}}
\put(13.00,40.00){\line(1,0){24.00}}
\put(15.00,37.50){\line(1,0){5.00}}
\put(15.00,42.50){\line(1,0){5.00}}
\put(20.00,2.50){\line(1,0){10.00}}
\put(20.00,47.50){\line(1,0){10.00}}
\put(30.00,7.50){\line(1,0){5.00}}
\put(30.00,12.50){\line(1,0){5.00}}
\put(30.00,37.50){\line(1,0){5.00}}
\put(30.00,42.50){\line(1,0){5.00}}
\put(35.00,17.50){\line(1,0){10.00}}
\put(35.00,32.50){\line(1,0){10.00}}
\put(45.00,22.50){\line(1,0){5.00}}
\put(45.00,27.50){\line(1,0){5.00}}

\put(25.00,-2.00){\line(0,1){54.00}}
\put(22.50,0.00){\line(0,1){5.00}}
\put(27.50,0.00){\line(0,1){5.00}}
\put(17.50,5.00){\line(0,1){10.00}}
\put(32.50,5.00){\line(0,1){10.00}}
\put(10.00,13.00){\line(0,1){24.00}}
\put(7.50,15.00){\line(0,1){5.00}}
\put(12.50,15.00){\line(0,1){5.00}}
\put(40.00,13.00){\line(0,1){24.00}}
\put(37.50,15.00){\line(0,1){5.00}}
\put(42.50,15.00){\line(0,1){5.00}}
\put(2.50,20.00){\line(0,1){10.00}}
\put(47.50,20.00){\line(0,1){10.00}}
\put(7.50,30.00){\line(0,1){5.00}}
\put(12.50,30.00){\line(0,1){5.00}}
\put(37.50,30.00){\line(0,1){5.00}}
\put(42.50,30.00){\line(0,1){5.00}}
\put(17.50,35.00){\line(0,1){10.00}}
\put(32.50,35.00){\line(0,1){10.00}}
\put(22.50,45.00){\line(0,1){5.00}}
\put(27.50,45.00){\line(0,1){5.00}}

\put(2.50,22.50){\circle*{1.00}}
\put(2.50,25.00){\circle*{1.00}}
\put(2.50,27.50){\circle*{1.00}}
\put(7.50,17.50){\circle*{1.00}}
\put(7.50,32.50){\circle*{1.00}}
\put(10.00,17.50){\circle*{1.00}}
\put(10.00,25.00){\circle*{1.00}}
\put(10.00,32.50){\circle*{1.00}}

\put(12.50,17.50){\circle*{1.00}}
\put(12.50,32.50){\circle*{1.00}}
\put(17.50,7.50){\circle*{1.00}}
\put(17.50,10.00){\circle*{1.00}}
\put(17.50,12.50){\circle*{1.00}}
\put(17.50,37.50){\circle*{1.00}}
\put(17.50,40.00){\circle*{1.00}}

\put(17.50,42.50){\circle*{1.00}}
\put(22.50,2.50){\circle*{1.00}}
\put(22.50,47.50){\circle*{1.00}}
\put(25.00,2.50){\circle*{1.00}}
\put(25.00,10.00){\circle*{1.00}}
\put(25.00,25.00){\circle*{1.00}}
\put(25.00,40.00){\circle*{1.00}}

\put(25.00,47.50){\circle*{1.00}}
\put(27.50,2.50){\circle*{1.00}}
\put(27.50,47.50){\circle*{1.00}}
\put(32.50,7.50){\circle*{1.00}}
\put(32.50,10.00){\circle*{1.00}}
\put(32.50,12.50){\circle*{1.00}}
\put(32.50,37.50){\circle*{1.00}}
\put(32.50,40.00){\circle*{1.00}}
\put(32.50,42.50){\circle*{1.00}}
\put(37.50,17.50){\circle*{1.00}}
\put(37.50,32.50){\circle*{1.00}}
\put(40.00,17.50){\circle*{1.00}}
\put(40.00,25.00){\circle*{1.00}}
\put(40.00,32.50){\circle*{1.00}}

\put(42.50,17.50){\circle*{1.00}}
\put(42.50,17.50){\circle*{1.00}}
\put(42.50,32.50){\circle*{1.00}}
\put(47.50,22.50){\circle*{1.00}}
\put(47.50,25.00){\circle*{1.00}}
\put(47.50,27.50){\circle*{1.00}}

\put(26.00,26.00){\footnotesize $e$}

\put(73.00,40.00){\footnotesize ${\cal G}_{2}\cong {\mathbb H}_2$}
\put(72.00,3.00){\footnotesize ${\cal G}_{1}\cong {\mathbb H}_2$}

\put(70.00,-4.00){\footnotesize $x_{-2}$}
\put(76.00,-4.00){\footnotesize  $x_{-1}$}
\put(82.00,-4.00){\footnotesize $e_{1}$}
\put(87.00,-4.00){\footnotesize  $x_1$}
\put(93.00,-4.00){\footnotesize $x_2$}

\put(64.00,28.00){\footnotesize $y_{-2}$}
\put(64.00,33.00){\footnotesize  $y_{-1}$}
\put(64.00,38.00){\footnotesize $e_{2}$}
\put(64.00,43.00){\footnotesize  $y_1$}
\put(64.00,48.00){\footnotesize $y_2$}

\put(70.00,25.00){\line(0,1){26.00}}
\put(70.00,0.00){\line(1,0){26.00}}

\put(70.00,28.00){\circle*{1.00}}
\put(70.00,33.00){\circle*{1.00}}
\put(70.00,38.00){\circle*{1.00}}
\put(70.00,43.00){\circle*{1.00}}
\put(70.00,48.00){\circle*{1.00}}

\put(73.00,0.00){\circle*{1.00}}
\put(78.00,0.00){\circle*{1.00}}
\put(83.00,0.00){\circle*{1.00}}
\put(88.00,0.00){\circle*{1.00}}
\put(93.00,0.00){\circle*{1.00}}

\end{picture}\\
\begin{center}
{\it Fig. 4.}
Homogenous tree ${\mathbb H}_4\cong {\mathbb H}_2*{\mathbb H}_2$.\\[5pt]
\end{center}
In particular, this implies that
${\mathbb H}_4\cong {\mathbb H}_2*{\mathbb H}_2$, where ${\mathbb H}_2\cong {\mathbb Z}$ is
the homogenous tree of order 2, isomorphic  the two-way infinite path.
Also, ${\mathbb H}_2\cong {\mathbb Z}_{2}*{\mathbb Z}_{2}$ (see Fig.4)\\[10pt]
\myownsection
\begin{center}
{\sc 6. Free decomposition}
\end{center}
In this section we show that the adjacency matrix of the
free product of a finite family of rooted graphs is the sum of free independent copies
of the adjacency matrices of the factors of the product. Essentially, this fact is 
a natural consequence of free probability and one just needs to adapt the proof of Voiculescu,
which we do for the reader's convenience.

In order to give the explicit form of the adjacency matrix of the free product of graphs,
let us consider the Hilbert space ${\cal H}_{I}=l_{2}(*_{i\in I}V_{i})$ spanned by
$\delta(e)$ and vectors of the form
$$
\delta(w)\;\; {\rm for}\;\;
w=v_{1}v_{2}\ldots v_{n}\in *_{i\in I}V_{i}
$$
and let $\varphi$ be the vacuum expectation on $B({\cal H})$ given by
$\varphi (T)=\langle T\delta(e) , \delta(e)\rangle$. We have
$$
({\cal H}_{I}, \delta(e))\cong  *_{i\in I}
(l_{2}(V_{i}),\delta(e_i))
$$
where the RHS is understood as the free product of Hilbert spaces with
distinguished unit vectors [33].

In the sequel we will need the following subsets of $\;*_{i\in I}V_{i}$:
$$
W_{j}(n)=\{v_{1}v_{2}\ldots v_{n}\in *_{i\in I}V_{i}:\;\;v_{1}\notin V_{j}^{0}\}
$$
where $n\in {\mathbb N}$. Thus, $W_{j}(n)$ is the subset
consisting of words of lenght $n$ which do not begin with a letter from $V_{j}^{0}$.
We set
$$
W_{j}=\bigcup_{n=0}^{\infty}W_{j}(n)
$$
with $W_{j}(0)=\{e\}$ for every $j$.\\
\indent{\par}
{\sc Definition 6.1.}
Let $A_{i}$ denote the adjacency matrix of the rooted graph $({\cal G}_{i},e_{i})$,
where $i \in I$. Let us define their copies in
$*_{i\in I}({\cal G}_{i},e_i)$ by the formulas
$$
(A_{i}(n))_{w,w'}=
\left\{
\begin{array}{ll}
1 & {\rm if}\; \{w,w'\}=\{xu,x'u\}\;{\rm for}\;\{x,x'\}\in E_{j}\;{\rm and}\; u\in W_{j}(n-1)\\
0 & {\rm otherwise}
\end{array}
\right.
$$
where $n\in {\mathbb N}$. By $P_{i}(n)$ we denote
the canonical projection of ${\cal H}_{I}$
onto $l_{2}(W_{i}(n))$ for $n\geq 1$ with $P_{i}(0)$ denoting the
projection onto $l_{2}(e)={\mathbb C}\delta(e)$ for every $i\in I$.\\
\indent{\par}
{\sc Theorem 6.2.}
{\it  The adjacency matrix $A(*_{i\in I}{\cal G}_{i})$ of the free
product of graphs admits a decomposition of the form
$A(*_{i\in I}{\cal G}_{i})=\sum_{i\in I}A^{(i)}$, where}
\begin{equation}
A^{(i)}=\sum_{n=1}^{\infty}A_{i}(n)=\sum_{n=1}^{\infty}A_{i}P_{i}(n-1)
\end{equation}
{\it are free w.r.t. the vacuum expectation $\varphi$ and the action of $A_{i}$ in the second sum
is given by
$A_{i}\delta(xu)=\delta(x'u)$ whenever $\{x,x'\}\in E_{i}\;(i\in I)$.
Moreover, the series is strongly convergent for every $i\in I$.}\\
\indent{\par}
{\it Proof.}
First, let us observe that using local finitness
of ${\cal G}_{i}$ we can write
$$
A^{(i)}\delta(w)=
\sum_{
\stackrel{w'=x'u} {\scriptscriptstyle \{x,x'\}\in E_{i}}
}
\delta(w')
=A_{i}(n)\delta(w)
$$
whenever $w=xu$, where $u\in W_{i}(n-1)$
and we allow $x,x'$ to be arbitrary vertices from $V_{i}$,
thus if $x=e_1$, we have $e_1 u\equiv u$.
Moreover, observe that $A_{i}(m)\delta(w)=0$ for
such $w$ for any $m\neq n$.
Writing
$$
\delta(w)=
\left\{
\begin{array}{ll}
\delta(x)\otimes \delta(u) & {\rm whenever}\; w=xu\; {\rm and} \;x\in V_{i}^{0}\\
\delta(u) & {\rm if} \;w=e_iu\;
\end{array}
\right.
$$
we obtain
\begin{eqnarray*}
A^{(i)}\delta(w)&=&
\sum_{
\stackrel{w'=x'u} {\scriptscriptstyle \{x,x'\}\in E_{i},x'\neq e_i}}
\delta(x')\otimes \delta(u) + \1_{\{\{x,e_1\}\in E_{i}\}}\delta(u)\\
&=&
(A_{i}\delta(x))^{0}\otimes \delta(u) + \langle A_{i}\delta(x), \delta(e_i))\rangle \delta(u)
\end{eqnarray*}
if $w=xu$, $u\in W_{i}$, where $\1_{\{z\in A\}}=1$ if and only if $z\in A$ and
otherwise is zero.

We can observe now that $A^{(i)}=\lambda(A_{i})$, where $\lambda$ denotes the free product
representation of the free product ${\mathbb C}[A_{1}]*{\mathbb C}[A_{2}]$ on
$l_{2}(*_{i=1}^{n}V_{i})$ in the sense of Avitzour [3] and Voiculescu [33].
Therefore, the $A^{(i)}$ are free with respect to
$\varphi$. \hfill $\blacksquare$\\
\indent{\par}
As a consequence of the decomposition theorem, one can use the free additive convolutions [34]
to compute spectral distributions of free products of rooted graphs in terms of
spectral distributions of the factors. In particular, one obtains asymptotic spectral 
properties of free powers $({\cal G}, e)^{*n}$.\\
\indent{\par}
{\sc Corollary 6.3.}
{\it Let $A$ be the adjacency matrix of $({\cal G}, e)$ and
let $A^{*n}$ denote the adjacency matrix of $({\cal G}, e)^{*n}$. Then}
\begin{equation}
\lim_{n\rightarrow \infty}\varphi \left(\left(\frac{A^{*n}}{\sqrt{nk(e})}\right)^{2m}\right)=
c_{m}
\end{equation}
{\it where $c_{m}$ is the $m$-th Catalan number for $m\in {\mathbb N}$,
$c_{0}=1$ and $k(e)$ is the degree of the root $e$. The odd moments vanish.}\\
\indent{\par}
{\it Proof.}
Applying the decomposition of Theorem 6.2 into the sum of freely independent random variables, we have
$$
\varphi(A^{(i)})=\varphi(A_{i}(1))=0
$$
since ${\cal G}$ has no loops. Moreover,
$$
\varphi((A^{(i)})^{2})= \varphi ((A_{i}(1))^{2})=
\sum_{\stackrel{\{x,e\}\in E}{\scriptscriptstyle \{x',e\}\in E}}
\langle \delta(x), \delta(x')\rangle =k(e)
$$
where $E=E({\cal G})$.
Using the central limit theorem for free random  variables, we obtain
the moments of the Wigner measure, which are Catalan numbers for even moments
(odd moments vanish).
\hfill $\blacksquare$\\
\indent{\par}
{\sc Corollary 6.4.}
{\it Let $A_{i}$ be the adjacency matrix of $({\cal G}_{i},e_{i})$, $i\in I=\{1 ,\ldots , n\}$,
and let $\mu_{i}$ denote its spectral distribution, where $1\leq i \leq n$.
Then the spectral distribution of $A(*_{i=1}^{n}{\cal G}_{i})$ in the state
$\varphi$ is given by $\mu=\mu_{1} \boxplus \mu_{2}\boxplus \ldots \boxplus \mu_{n}$}.\\
\indent{\par}
{\it Proof.}
This is a direct consequence of the fact that
$A^{(1)},A^{(2)}, \ldots ,A^{(n)}$ are free w.r.t. the state $\varphi$.
\hfill $\blacksquare$\\
\indent{\par}
{\it Remark.}
The correspondence between the free product of graphs and free probability
(in particular, Corollary 6.4) can be applied to establish a connection between free
products of graphs and free additive convolutions of their spectral distributions.\\
\indent{\par}
{\sc Example 6.1.}
Using Corollary 6.4, we can find the spectral distributions $\nu_{n}$ of ${\mathbb T}_n$
associated with the root, 
where $n\geq 1$. Denote $\nu=\nu_{1}$. From Example 5.3 we obtain
$\nu_{n}=\nu\boxplus \nu\boxplus \ldots \boxplus \nu$ ($n$ times).
Knowing $G_{\nu}(z)=1/2(z-\sqrt{z^2-4})$ we can compute $G_{\nu}^{-1}(z)=z^{-1}+z$ and
thus $R_{\nu}(z)=z$, which gives $R_{\nu_{n}}(z)=nR_{\nu}(z)=nz$ by (3.7),
and this leads to $G_{\nu_{n}}^{-1}(z)=z^{-1}+nz$. Thus
\begin{equation}
G_{\nu_{n}}(z)=\frac{z-\sqrt{z^{2}-4n}}{2n}
\end{equation}
from which, by applying the Stieltjes inversion formula, 
we obtain the Wigner laws with densities
\begin{equation}
d\nu_{n}(x)=\frac{\sqrt{4n-x^{2}}}{2\pi n}dx
\end{equation}
with the support on $[-2\sqrt{n},2\sqrt{n}]$.\\
\indent{\par}
{\sc Example 6.2.}
In a similar manner, we can compute the spectral distribution $\mu_{n}$ of ${\mathbb H}_n$
associated with the root (or any vertex, since all vertices are equivalent).
Denote $\mu=\mu_{1}$. Using Example 5.4 and Corollary 6.4, we obtain
$\mu_{n}=\mu\boxplus \mu\boxplus \ldots \boxplus \mu$ ($n$ times) and thus
it can be computed using the R-transforms. Namely, from
$G_{\mu}(z)=z/(z^{2}-1)$, we get
$$
R_{\mu}(z)=\frac{-1+\sqrt{1+4z^{2}}}{2z}
$$
and thus $R_{\mu_{n}}=nR_{\mu}$ using (3.7).
Therefore, the Cauchy transforms of the measures $\mu_{n}$ are
\begin{equation}
G_{\mu_{n}}(z)=\frac{(2-n)z+n\sqrt{z^{2}-4(n-1)}}{2(z^{2}-n^{2})}
\end{equation}
which, with the help of the Stieltjes inversion formula, give the 
(absolutely continuous) measures with densities
\begin{equation}
d\mu_{n}(x)=\frac{n\sqrt{4(n-1)-x^{2}}}{2\pi (n^{2}-x^{2})}dx
\end{equation}
supported on $[-2\sqrt{n-1}, 2\sqrt{n-1}]$.
In particular, the spectral distribution of
${\mathbb H}_2\cong {\mathbb Z}$ in the vacuum state $\varphi$ associated with the vertex $0$ 
is the arcsine law $d\mu_{2}(x)=1/(\pi\sqrt{4-x^{2}})dx$.\\[10pt]

\myownsection
\begin{center}
{\sc 7. Orthogonal decompositions of branches}
\end{center}
Let us look at the concept of `branches' of the free product of graphs
introduced by Quenell [30]. They correspond to the so-called `subordination functions' studied 
first by Voiculescu [36] and Biane [4]. We rely on the recent general study of
the free additive convolution and its decompositions given in [21], where we refer
the reader for the main concepts, like s-freeness, as well as general proofs.  \\
\indent{\par}
{\sc Definition 7.1.}
Let $(V_{i},e_{i})_{i\in I}$ be a finite family of rooted sets. 
By the {\it branch of $*_{i\in I}(V_{i},e_{i})$ subordinate to $(V_{j},e_{j})$},
where $j\in I$, we shall understand the rooted set $(S_{j},e)$, where
$$
S_{j}=
\{e\} \cup \{v_{1}v_{2}\ldots v_{m}\in *_{i\in I}V_{i}:\;\;
v_{m}\in V_{j}^{0},\; m\in {\mathbb N}\}
$$
is the subset of $*_{i\in I}V_{i}$ consisting of the empty word and words which end with
a letter from $V_{j}^{0}$. \\
\indent{\par}
{\sc Definition 7.2.}
Let $({\cal G}_{i},e_{i})_{i\in I}$ be a finite family of rooted graphs.
By the {\it branch of $*_{i\in I}({\cal G}_{i}, e_{i})$ subordinate to} 
$({\cal G}_{j},e_{j})$, where $j\in I$,
we shall understand the rooted graph $({\cal B}_{j},e)$,
where ${\cal B}_{j}\equiv{\cal B}_{j}(({\cal G}_{i})_{i\in I})$
is the subgraph of $*_{i\in I}{\cal G}_{i}$ restricted to the set
$S_{j}$ defined above. As before, we often omit the roots in the notations.\\

\unitlength=1mm
\special{em.linewidth 0.5pt}
\linethickness{0.5pt}
\begin{picture}(140.00,60.00)(-15.00,5.00)

\put(10.00,10.00){\line(1,2){5.00}}
\put(20.00,10.00){\line(-1,2){5.00}}
\put(30.00,10.00){\line(1,2){5.00}}
\put(40.00,10.00){\line(-1,2){5.00}}

\put(50.00,10.00){\line(1,2){5.00}}
\put(60.00,10.00){\line(-1,2){5.00}}
\put(70.00,10.00){\line(1,2){5.00}}
\put(80.00,10.00){\line(-1,2){5.00}}

\put(15.00,20.00){\line(2,3){10.00}}
\put(35.00,20.00){\line(-2,3){10.00}}
\put(55.00,20.00){\line(2,3){10.00}}
\put(75.00,20.00){\line(-2,3){10.00}}

\put(25.00,35.00){\line(4,3){20.00}}
\put(65.00,35.00){\line(-4,3){20.00}}

\put(10.00,10.00){\circle*{1.50}}
\put(20.00,10.00){\circle*{1.50}}
\put(30.00,10.00){\circle*{1.50}}
\put(40.00,10.00){\circle*{1.50}}
\put(50.00,10.00){\circle*{1.50}}
\put(60.00,10.00){\circle*{1.50}}
\put(70.00,10.00){\circle*{1.50}}
\put(80.00,10.00){\circle*{1.50}}

\put(10.00,10.00){\line(-1,-2){2.00}}
\put(10.00,10.00){\line(1,-2){2.00}}
\put(20.00,10.00){\line(-1,-2){2.00}}
\put(20.00,10.00){\line(1,-2){2.00}}
\put(30.00,10.00){\line(-1,-2){2.00}}
\put(30.00,10.00){\line(1,-2){2.00}}
\put(40.00,10.00){\line(-1,-2){2.00}}
\put(40.00,10.00){\line(1,-2){2.00}}

\put(50.00,10.00){\line(-1,-2){2.00}}
\put(50.00,10.00){\line(1,-2){2.00}}
\put(60.00,10.00){\line(-1,-2){2.00}}
\put(60.00,10.00){\line(1,-2){2.00}}
\put(70.00,10.00){\line(-1,-2){2.00}}
\put(70.00,10.00){\line(1,-2){2.00}}
\put(80.00,10.00){\line(-1,-2){2.00}}
\put(80.00,10.00){\line(1,-2){2.00}}

\put(15.00,20.00){\circle*{1.50}}
\put(35.00,20.00){\circle*{1.50}}
\put(55.00,20.00){\circle*{1.50}}
\put(75.00,20.00){\circle*{1.50}}

\put(25.00,35.00){\circle*{1.50}}
\put(65.00,35.00){\circle*{1.50}}

\put(45.00,50.00){\circle*{1.50}}

\put(4.00,10.00){\footnotesize $xyx$}
\put(13.00,10.00){\footnotesize $x'yx$}
\put(22.50,10.00){\footnotesize $xy'x$}
\put(32.00,10.00){\footnotesize $x'y'x$}
\put(51.00,10.00){\footnotesize $xyx'$}
\put(61.00,10.00){\footnotesize $x'yx'$}
\put(72.00,10.00){\footnotesize $xy'x'$}
\put(81.00,10.00){\footnotesize $x'y'x'$}

\put(10.00,20.00){\footnotesize $yx$}
\put(37.00,20.00){\footnotesize $y'x$}
\put(49.00,20.00){\footnotesize $yx'$}
\put(77.00,20.00){\footnotesize $y'x'$}

\put(20.00,35.00){\footnotesize $x$}
\put(67.00,35.00){\footnotesize $x'$}
\put(44.00,52.00){\footnotesize $e$}

\put(100.00,10.00){\line(1,1){6.00}}
\put(112.00,10.00){\line(-1,1){6.00}}

\put(100.00,10.00){\circle*{1.50}}
\put(112.00,10.00){\circle*{1.50}}
\put(106.00,16.00){\circle*{1.50}}

\put(99.00,6.00){\footnotesize $x$}
\put(105.00,18.00){\footnotesize $e_{1}$}
\put(111.00,6.00){\footnotesize $x'$}

\put(100.00,30.00){\line(1,1){6.00}}
\put(112.00,30.00){\line(-1,1){6.00}}

\put(100.00,30.00){\circle*{1.50}}
\put(112.00,30.00){\circle*{1.50}}
\put(106.00,36.00){\circle*{1.50}}

\put(99.00,26.00){\footnotesize $x$}
\put(105.00,38.00){\footnotesize $e_{1}$}
\put(111.00,26.00){\footnotesize $x'$}

\put(96.00,15.00){\footnotesize ${\cal G}_{1}$}
\put(96.00,35.00){\footnotesize ${\cal G}_{2}$}

\end{picture}\\
\begin{center}
{\it Fig. 5.}
Binary tree ${\mathbb T}_2\cong {\cal B}_{1}({\cal G}_{1}*{\cal G}_{2})$.\\[5pt]
\end{center}
\indent{\par}
In the case of two graphs, it is easy to see that
${\cal G}_{1}*{\cal G}_{2}$ consists of two branches ${\cal B}_{1}$
and ${\cal B}_{2}$ with common root $e$.
The branch ${\cal B}_{1}={\cal B}_{1}({\cal G}_{1},{\cal G}_{2})$
`begins' with a copy of ${\cal G}_{1}$ and the branch ${\cal B}_{2}={\cal B}_{2}({\cal G}_{1},{\cal G}_{2})$ `begins' with a copy of ${\cal G}_{2}$. 
For instance, in the case of a binary tree ${\mathbb T}_2$ (Fig.3), 
the branches ${\cal B}_{1}$ and ${\cal B}_{2}$ are the left and right
`halves' of ${\mathbb T}_2$, respectively. In the case of ${\mathbb H}_4$ (Fig.4), the branch ${\cal B}_{1}$
(${\cal B}_{2}$) of ${\mathbb H}_4$ consists of the horizontal (vertical) `diagonal' together with
all `leaves' attached to all its vertices.
However, the $n$-ary tree can itself be viewed as a branch of another free product. Moreover,
it is then constructed in a `distance-adapted' manner, i.e. natural truncations
of the product lead to natural truncations of the tree (see Example 7.1 and Fig.5).\\
\indent{\par}
{\sc Example 7.1.}
Take the same two graphs as in Example 5.2, but choose the roots $e_{1}$ and
$e_{2}$ in such a way that $x\sim e_{1}\sim x'$ and $y\sim e_{2}\sim y'$
and let ${\cal B}_{1}$ and ${\cal B}_{2}$ be the branches of the free product 
${\cal G}_{1}*{\cal G}_{2}$ subordinate to ${\cal G}_{1}$
and ${\cal G}_{2}$, respectively. 
Figure 5 shows that ${\mathbb T}_2\cong {\cal B}_{1}(\cong {\cal B}_{2})$.\\
\indent{\par}
In a similar way one can obtain the $n$-ary rooted tree as a branch of a free
product of two copies of the `fork' graph with $n+1$ vertices
$e,x_{1},\ldots , x_{n}$ (i.e. such such that $e\sim x_{k}$ for all $1\leq k \leq n$).\\
\indent{\par}
Motivated by [21], we can view the branches of ${\cal G}_{1}*{\cal G}_{2}$ as 
products of rooted graphs. The needed notion of a product corresponds to 
{\it freeness with subordination}, or {\it s-freeness}, introduced and studied there. 
Moreover, their adjacency matrices, $A({\cal B}_{1})$ and $A({\cal B}_{2})$, 
can be decomposed as the sum of components which are `free with subordination' 
(or, simply `s-free') w.r.t. a pair of states $(\varphi, \psi)$. Note that 
freeness with subordination is quite different from conditional freeness [5], 
although it also involves two functionals (states).
\\
\indent{\par}
{\sc Definition 7.3.}
Let $({\cal A},\varphi, \psi)$ be a unital algebra with a pair
of linear normalized functionals.
Let ${\cal A}_{1}$ be a unital subalgebra of ${\cal A}$ and let
${\cal A}_{2}$ be a non-unital subalgebra with an `internal' unit $1_{2}$, i.e.
$1_{2}b=b=b1_{2}$ for every $b\in {\cal A}_{2}$.
We say that the pair $({\cal A}_{1},{\cal A}_{2})$ is {\it free with subordination}, or 
simply {\it s-free},
with respect to $(\varphi, \psi)$ if $\psi(1_{2})=1$ and it holds that\\
\indent{\par}
(i)
$\varphi(a_{1}a_{2}\ldots a_{n})=0$ whenever $a_{j}\in {\cal A}_{i_{j}}^{0}$ and
$i_{1}\ne i_{2}\neq \ldots \neq i_{n}$
\indent{\par}
(ii)
$\varphi(w_{1}1_{2}w_{2})=\varphi(w_{1}w_{2})-\varphi(w_{2})\varphi(w_{2})$
for any $w_{1},w_{2}\in {\rm alg}({\cal A}_{1}, {\cal A}_{2})$,\\[5pt]
where ${\cal A}_{1}^{0}={\cal A}_{1}\cap {\rm ker}\varphi$ and ${\cal A}_{2}^{0}={\cal A}_{2}\cap {\rm ker}\psi$.
We say that the pair $(a,b)$ of random variables from ${\cal A}$
is {\it s-free} with respect to $(\varphi, \psi)$
if the pair of algebras generated by these random variables is s-free with respect to $(\varphi, \psi)$.\\
\indent{\par}
The notion of s-freeness resembles freeness - in the GNS representation, the corresponding
product of Hilbert spaces is spanned by the unit (vacuum) vector $\xi$ and
simple tensors of the form ${\cal H}_{i_{1}}\otimes {\cal H}_{i_{2}}\otimes \ldots \otimes {\cal H}_{i_{n}}$,
where $i_{1}\neq i_{2}\neq \ldots \neq i_{n}=1$. The branches, which in this context 
replace free products of graphs, can also be decomposed along the lines of Theorem 6.2
and can be called `s-free products' of ${\cal G}_{1}$ and ${\cal G}_{2}$, or vice versa.
Using a similar notation, we obtain the following decomposition theorem (cf. [21]).\\
\indent{\par}
{\sc Theorem 7.1.}
{\it The adjacency matrix of the branch ${\cal B}_{1}\equiv{\cal B}_{1}({\cal G}_{1},{\cal G}_{2})$
can be decomposed as the sum $A({\cal B}_{1})=A^{(1)}+A^{(2)}$, where the strongly convergent series}
\begin{equation}
A^{(1)}=\sum_{n\; {\rm odd}}A_{1}(n), \;\;\;
A^{(2)}=\sum_{n\; {\rm even}}A_{2}(n),
\end{equation}
{\it are s-free w.r.t.
$(\varphi, \psi)$, where $\varphi(.)=\langle .\delta(e),\delta(e)\rangle$ and 
$\psi(.)=\langle . \delta(v),\delta(v)\rangle$ and $v\in V_{1}^{0}$.
An analogous decomposition holds for the branch ${\cal B}_{2}({\cal G}_{1}*{\cal G}_{2})$
with the summations over odd and even $n$ interchanged.} \\
\indent{\par}
{\it Proof.}
We refer the reader to [21], where it was shown, in a general Hilbert space setting, 
that sums of operators of the above type are s-free 
w.r.t. $(\varphi, \psi)$ (one has to verify conditions (i)-(ii) of Definition 7.3, 
and in the case of graphs, it is basically the same proof).
\hfill $\blacksquare$\\
\indent{\par}
In order to `decompose completely' the branches, by which we mean to decompose them in terms
of graphs ${\cal G}_{1}$ and ${\cal G}_{2}$, we will interpret (7.1) in terms
of an inductive limit of a sequence of graphs which resembles (but is not the same as)
the sequence of $m$-free products approximating the free product. In this fashion we will obtain 
the `complete' {\it orthogonal decomposition} of branches given by the following theorem. \\
\indent{\par}
{\sc Theorem 7.2.}
{\it The branch ${\cal B}_{1}$ is the inductive limit of the sequence ${\cal G}_{1}\vdash_{m} {\cal G}_{2}$
given by the recursion}
$$
{\cal G}_{1}\vdash_{1}{\cal G}_{2}={\cal G}_{1}\vdash {\cal G}_{2}, \;\;\;\;
{\cal G}_{1}\vdash_{m}{\cal G}_{2}={\cal G}_{1}\vdash ({\cal G}_{2}\vdash_{m-1} {\cal G}_{2}),
$$
{\it where $m>1$. An analogous statement holds for the branch ${\cal B}_{2}$.}\\
\indent{\par}
{\it Proof.}
Without loss of generality, consider branch ${\cal B}_{1}$.
Our sequence of iterates will remind the inductive way of defining the 
free product of graphs given in [39], although it is asymmetric with respect to
${\cal G}_{1}$ and ${\cal G}_{2}$.
Recall that ${\cal B}_{i}$ `begins' with a copy of ${\cal G}_{i}$.
Therefore, let ${\cal B}_{1}(0)$ be equal to ${\cal G}_{1}$ and choose its root to be $e_{1}$.
To get ${\cal B}_{1}(1)$, to every vertex of ${\cal G}_{1}$ but the root we glue by its root 
a copy of ${\cal B}_{2}(0)$. In such a graph we again choose 
the root $e_{1}$. This gives a rooted graph $({\cal B}_{1}(1),e_{1})$, 
which is, in fact ${\cal G}_{1}\vdash {\cal G}_{2}$. In a similar fashion we obtain 
$({\cal B}_{2}(1),e_{2})$. Now, note that the $m$-th approximant of the branch ${\cal B}_{1}$ 
is obtained by glueing by its root a copy of $({\cal B}_{2}(m-1), e_{2})$ to every vertex of 
$({\cal G}_{1},e_{1})$ but the root. In other words, we obtain
$$
{\cal B}_{1}(m)={\cal G}_{1}\vdash {\cal B}_{2}(m-1)\;\;\;{\rm and} \;\;\;
{\cal B}_{2}(m)={\cal G}_{2}\vdash {\cal B}_{1}(m-1)
$$
for $m\geq 1$. It is clear that the inductive limits of our iterates 
give the branches, namely
$$
{\cal B}_{i}=\bigcup_{m \geq 0}{\cal B}_{i}(m)
$$
for $i=1,2$ (with the root $e_{i}$), and this proves the assertion. \hfill $\blacksquare$\\
\indent{\par}
In oder to obtain spectral distributions of the branches, one takes 
a sequence of alternating iterates of orthogonal convolutions -
this method was introduced in [21], but here, in the graph context, is
especially appealing and easy to justify.\\
\indent{\par}
{\sc Corollary 7.3.}
{\it If $\mu$ and $\nu$ are spectral distributions of ${\cal G}_{1}$ and ${\cal G}_{2}$, 
respectively, then the spectral distribution of ${\cal G}_{1}\vdash_{m} {\cal G}_{2}$
is given by $\mu \vdash_{m}\nu$, where the sequence $(\mu \vdash_{m}\nu)_{m\in {\mathbb N}}$
of distributions is given by the recursion
$$
\mu \vdash_1 \nu=\mu  \vdash  \nu, \;\;\;
\mu \vdash_{m} \nu= \mu\vdash (\nu \vdash_{m-1}\mu)
$$
where $m>1$. If ${\cal G}_{1}$ and ${\cal G}_{2}$ are uniformly locally finite,
the spectral distribution of the branch ${\cal B}_{1}$ is given by the weak limit
$\mu\boxright\nu := w-\lim_{m\rightarrow \infty} \mu\vdash_{m} \nu .$
An analogous statement holds for the branch ${\cal B}_{2}$.}\\
\indent{\par}
{\it Proof.}
The spectral distribution of ${\cal G}_{1}\vdash_{m} {\cal G}_{2}$ 
is given by $\mu \vdash_{m}\nu$ by Theorem 4.1. 
Now, observe that moments of the same order $k$, where $k\leq 2m$,
in all graphs ${\cal G}_{1}\vdash_n{\cal G}_{2}$ (computed w.r.t. the root $e$), with $n\geq m$, are 
equal. This is because in the orthogonal product of graphs no copy of the second graph 
is glued to the root of the first graph and therefore, the distance from the root $e$ in 
${\cal G}_{1}\vdash_{m}{\cal G}_{2}$ at which the graph differs from ${\cal G}_{1}\vdash_{m-1}{\cal G}_{2}$ 
is equal to $m+1$. Therefore, the sequence of moments of $\mu \vdash_{m}\nu$
converges to the moments of ${\cal B}_{1}$. If ${\cal G}_{1}$ and ${\cal G}_{2}$ are
uniformly locally finite, this implies weak convergence of measures.\hfill $\blacksquare$\\
\indent{\par}
{\sc Corollary 7.4.}
{\it Under the assumptions of Corollary 7.3, the K-transform of $\mu \,\boxright \,\nu$
can be expressed as}
$$
K_{\mu\,\boxright\,\nu}(z)=K_{\mu}(z-K_{\nu}(z-K_{\mu}(z-K_{\nu}(\ldots ))))
$$
{\it where the right-hand side is understood as the uniform limit on compact subsets
of the complex upper half-plane. The K-transform of $\mu\vdash_{m}\nu$
is obtained by a truncation of the above formula to $m+1$ alternating transforms.}\\
\indent{\par}
{\it Proof.}
Since weak convergence of measures implies uniform convergence of the Cauchy transform
on compact subsets of the complex upper half-plane, the assertion follows from
a repeated application of (3.9) and Corollary 7.3. \hfill $\blacksquare$\\
\indent{\par}
Let us note that the above `continued composition formula' is very convenient for
computing the K-,F-, or G-transforms of $\mu\,\boxright \,\nu$. Essentially,
for all examples of graphs whose free products have been studied so far, it provides
a tool which immediately gives the continued fractions of their s-free product. It also explains
why 2-periodic and mixed periodic Jacobi continued fractions are so typical
in the context of free products.\\
\indent{\par}
{\sc Example 7.2.}
Consider two rooted graphs ${\cal G}_{1}$, ${\cal G}_{2}$, whose spectral distributions
$\mu, \nu$ are associated with reciprocal Cauchy transforms of the form
$$
F_{\mu}(z)=z-\alpha_{0}-\frac{\omega_{0}}{z-\alpha_{1}}, \;\;\;
F_{\nu}(z)=z-\beta_{0}-\frac{\gamma_{0}}{z-\beta_{1}},
$$
respectively (this includes ${\mathbb K}_n$ and ${\mathbb F}_n$, whose free
products were studied by other authors and also in Section 11). 
From Corollary 7.4 we easily obtain the K-transform
$$
K_{\mu\,\boxright \,\nu}(z)=
\alpha_{0}+\cfrac{\omega_{0}}{z-\alpha_{1}-\beta_{0}-\cfrac{\gamma_{0}}{z-\alpha_{0}-\beta_{1}-\cfrac{\omega_{0}}{z-\alpha_{1}-\beta_{0}-\cfrac{\gamma_{0}}{\ldots}}}}
$$
and thus, in view of (3.3), the distribution $\mu \,\boxright \,\nu$ of branch ${\cal B}_{1}$
is associated with the sequences of Jacobi parameters
$$
\alpha=(\alpha_{0},\alpha_{1}+\beta_{0}, \alpha_{0}+\beta_{1}, \alpha_{1}+\beta_{0}, \ldots), \;\;\;
\omega=(\omega_{0},\gamma_{0},\omega_{0},\gamma_{0}, \ldots)
$$
which correspond to the so-called mixed periodic Jacobi continued fraction [14].
For details on the corresponding measures, see [14]. In particular, if 
${\cal G}_{1}={\cal G}_{2}={\mathbb K}_2$ (3-vertex complete graph), we have
$\alpha_{0}=\beta_{0}=0$, $\alpha_{1}=\beta_{1}=1$ 
$\omega_{0}=\gamma_{0}=2$ which gives $\mu\,\boxright \,\nu$ associated
with the sequences of Jacobi parameters $\alpha=(0,1,1,\ldots)$ and $\omega=(2,2,\ldots)$.
Its Cauchy transform is 
$$
G(z)=\frac{z+1-\sqrt{z^{2}-2z-7}}{2z+4}
$$
and the measure has density $d\mu(x)=\sqrt{7+2x-x^{2}}/(\pi(2x+4))$ on 
the interval $[1-2\sqrt{2},1+2\sqrt{2}]$.\\[10pt]

\myownsection
\begin{center}
{\sc 8. `Complete' decompositions of free products}
\end{center}
In this section we derive new decompositions of the free product of graphs,
which are based on the orthogonal decomposition of branches of Section 7.
We rely on the general theory of the decompositions of the
free additive convolution [21] and apply it to the context of graph products.

We start from two lemmas, which rephrase the results of Quenell
using the language of quantum probability. This reduces certain
proofs presented in [30] to the basic properties of 
the monotone and boolean convolutions.\\
\indent{\par}
{\sc Lemma 8.1.}
{\it The free product of rooted graphs admits the decomposition}
\begin{equation}
{\cal G}_{1}*{\cal G}_{2}
\cong
{\cal B}_{1}\star {\cal B}_{2}
\end{equation}
{\it which we call the star decomposition of ${\cal G}_{1}*{\cal G}_{2}$}.\\
\indent{\par}
{\it Proof.}
Notice that $V_{1}*V_{2}=S_{1}\cup S_{2}$.
Moreover, if follows immediately from Definition 7.2 that 
the free product ${\cal G}_{1}*{\cal G}_{2}$ is obtained by glueing together the branches
${\cal B}_{1}$ and ${\cal B}_{2}$ at their roots. From the definition of the star
product we know that this is the star product of ${\cal B}_{1}$ and ${\cal B}_{2}$.
\hfill $\blacksquare$\\
\indent{\par}
{\sc Lemma 8.2.}
{\it The free product of rooted graphs admits the decompositions}
\begin{equation}
{\cal G}_{1}*{\cal G}_{2}
\cong {\cal G}_{1}   \vartriangleright {\cal B}_{2} \cong 
{\cal G}_{2}   \vartriangleright {\cal B}_{1} 
\end{equation}
{\it which we call the comb decompositions of ${\cal G}_{1}*{\cal G}_{2}$.}\\
\indent{\par}
{\it Proof.}
Without loss of generality, consider the first relation.
Observe that we can view the branch ${\cal B}_{1}$ as one replica of graph ${\cal G}_{1}$, 
to which we glue `orthogonally' replicas of branch ${\cal B}_{2}$. Therefore
$$
{\cal B}_{1}\cong {\cal G}_{1}\vdash {\cal B}_{2} \;\;{\rm and}\;\;
{\cal B}_{2}\cong {\cal G}_{2}\vdash {\cal B}_{1}.
$$
By Theorem 8.1, we can obtain ${\cal G}_{1}*{\cal G}_{2}$ by glueing ${\cal B}_{1}$
and ${\cal B}_{2}$ together at their roots identified with $e$. 
Equivalently, (one replica of) ${\cal B}_{2}$ is glued to $e$ and branch ${\cal B}_{1}$  
is replaced by ${\cal G}_{1}\vdash {\cal B}_{2}$, which means that a replica of 
${\cal B}_{2}$ is glued to every vertex $v\in V_{1}^{0}\subset V_{1}*V_{2}$.
In other words, a replica of ${\cal B}_{2}$ is glued to every vertex of
$V_{1}^{0}\cup \{e\}\cong V_{1}$, which gives the comb product of ${\cal G}_{1}$
and ${\cal B}_{2}$, which proves our assertion. \hfill $\blacksquare$\\
\indent{\par}
{\sc Corollary 8.3.}
{\it The following relations hold:}
\begin{eqnarray*}
F_{\mu \,\boxplus \,\nu}(z)
&=& F_{\mu}(F_{\nu \,\boxright \,\mu}(z))+F_{\nu}(F_{\mu\,\boxright \,\nu}(z))-z\\
F_{\mu\,\boxplus\, \nu}(z)
&=&
F_{\mu}(F_{\nu \,\boxright\, \mu}(z))
=
F_{\nu}(F_{\mu \,\boxright \,\nu}(z))
\end{eqnarray*}
{\it where $\mu$ and $\nu$ are spectral distributions of rooted graphs
${\cal G}_{1}$ and ${\cal G}_{2}$.}\\
\indent{\par}
{\it Proof.}
These are straightforward consequences of Theorems 8.1-8.2 and formulas given by
(3.10) and (3.8).\hfill $\blacksquare$\\
\indent{\par}
{\sc Remark 8.1.}
In terms of moment generating functions $M_{\mu}(z)$, related to $F_{\mu}(z)$ by
(3.2) and (3.4), formulas of Corollary 8.3 give the results of Quenell [30] (see also
[24]). In our notation, $M_{\mu}(z)$ correspond to return generating functions of type $R_{e}(z)$ 
(from which the first return generating functions of type $S_{e}(z)$ are easily obtained).
Moreover, these results can be easily generalized to a finite number of rooted graphs.\\
\indent{\par}
Let us now use the `complete' orthogonal decomposition of branches (Theorem 7.2) 
and Lemmas 8.1-8.2 to derive `complete' decompositions of the free product of graphs. We begin with 
a decomposition of $m$-free products.\\
\indent{\par}
{\sc Theorem 8.4.}
{\it Let ${\cal G}_{1}$ and ${\cal G}_{2}$ be rooted graphs with spectral distributions
$\mu$ and $\nu$. Then their $m$-free product admits the decomposition}
$$
{\cal G}_{1}*^{m}{\cal G}_{2} 
=
({\cal G}_{1}\vdash_{m}{\cal G}_{2})\star ({\cal G}_{2}\vdash_{m} {\cal G}_{1})
$$
{\it called the star-orthogonal decomposition, and its spectral distribution is given by 
the m-free convolution $\mu \boxplus_{m} \nu :=(\mu \vdash_{m}\nu) \uplus (\nu\vdash_{m}\mu)$.}\\
\indent{\par}
{\it Proof.}
The proof consists in describing how to obtain the iterates of the free product
${\cal G}_{1}*{\cal G}_{2}$ in an inductive manner by appropriate glueing. 
Thus, in the first step we obtain
${\cal G}_{1}*^{1} {\cal G}_{2}$ by glueing one copy of ${\cal G}_{1}$
to one copy of ${\cal G}_{2}$ by means of identifying 
$e_1$ with $e_2$ and choosing it to be the root $e$. This gives
$({\cal G}_{1}\vdash {\cal G}_{2}) \star ({\cal G}_{2}\vdash {\cal G}_{1})$.
Note that in the $m$-th step we can obtain the graph ${\cal G}_{1}*^{m} {\cal G}_{2}$
by glueing a copy of ${\cal G}_{2}*^{m-1} {\cal G}_{1}$ 
to every vertex of ${\cal G}_{1}$ but the root $e_{1}$, and vice versa, 
a copy of ${\cal G}_{1}*^{m-1} {\cal G}_{2}$ to every vertex of 
${\cal G}_{2}$ but the root $e_{2}$, and then, by glueing the two graphs 
obtained in that way at their roots ($e_{1}$ and $e_{2}$, respectively).
These rules of glueing correspond to the orthogonal and star products and 
thus the first assertion is proved. The second assertion is then a consequence
of Corollary 7.2. \hfill $\blacksquare$\\
\indent{\par}
Below we state our results on the decompositions of the free product of rooted uniformly 
locally finite graphs.
The limits of products of rooted graphs are understood as inductive limits
(towers of graphs with the same root).
Results concerning convolutions have been proven in [21] for compactly supported
probability measures (see also [18], where it is shown that 
$\mu\,\boxplus \,\nu=w-\lim_{m\rightarrow \infty}\mu \boxplus_{m}\nu$, 
with a different, purely algebraic, definition of the $m$-free convolution). \\
\indent{\par}
{\sc Theorem 8.5.}
{\it Let ${\cal G}_{1}$ and ${\cal G}_{2}$ be rooted graphs with spectral distributions
$\mu$ and $\nu$. Then their free product admits the decomposition}
$$
{\cal G}_{1}*{\cal G}_{2}\cong ({\cal G}_{1}\vdash ({\cal G}_{2}\vdash ({\cal G}_{1}\vdash \ldots )))\star
({\cal G}_{2}\vdash ({\cal G}_{1}\vdash ({\cal G}_{2}\vdash \ldots )))
$$
{\it called the star-orthogonal decomposition. If ${\cal G}_{1}$ and ${\cal G}_{2}$ are
uniformly locally finite, its spectral distribution is given by 
$\;\mu \,\boxplus \,\nu= w-\lim_{m\rightarrow \infty} ((\mu \vdash_{m}\nu) \uplus (\nu\vdash_{m}\mu))$}.\\
\indent{\par}
{\it Proof.}
The first statement follows from Theorem 7.1 and Lemma 8.1. The weak limit formula for $\mu\, \boxplus \,\nu$
is a consequence of Corollary 7.3.\hfill $\blacksquare$\\
\indent{\par}
{\sc Theorem 8.6.}
{\it Under the assumptions of Theorem 8.5, the free product of rooted graphs admits the decomposition}
$$
{\cal G}_{1}*{\cal G}_{2}
\cong
{\cal G}_{1}   \vartriangleright ({\cal G}_{2}\vdash({\cal G}_{1}\vdash({\cal G}_{2}\vdash \ldots )))
$$
{\it called the comb-orthogonal decomposition. If ${\cal G}_{1}$ and ${\cal G}_{2}$
are uniformly locally finite,  its spectral distribution is
given by $\;\mu\,\boxplus\, \nu =w-\lim_{m\rightarrow \infty}
\mu \vartriangleright (\nu \vdash_{m}\mu)$.}\\
\indent{\par}
{\it Proof.}
The first statement follows from Theorem 7.2 and Lemma 8.2. The formula
for $\mu \,\boxplus \,\nu$ is a consequence of Corollary 7.3.\hfill $\blacksquare$\\
\indent{\par}
{\sc Corollary 8.7.}
{\it Under the assumptions of Corollary 7.4, the Cauchy transform of $\mu \,\boxplus \,\nu$
can be expressed as}
$$
G_{\mu\,\boxplus\,\nu}(z)=G_{\mu}(z-K_{\nu}(z-K_{\mu}(z-K_{\nu}(\ldots ))))
$$
{\it where the right-hand side is understood as the uniform limit on compact subsets
of the complex upper half-plane.}\\
\indent{\par}
{\it Proof.}
This is a consequence of repeated application of (3.9) and Theorem 8.6.\hfill $\blacksquare$\\
\indent{\par}
If ${\cal G}_{1}$ or ${\cal G}_{2}$ is not uniformly locally finite,
the statement of Theorem 8.6 concerning spectral distributions holds 
in the weaker sense of convergence of moments.\\
\indent{\par}
{\sc Example 8.1.}
Take two graphs of type given in Example 7.2. In contrast to the 
s-free product considered there, one cannot obtain an explicit formula for the continued 
Jacobi fraction corresponding to the free product (essentially, due to the presence of $G_{\mu}$
in the beginning of the above formula). Instead, we immediately obtain the algebraic formula
$$
G_{\mu\,\boxplus \,\nu}(z)=\cfrac{1}{z-\alpha_{0}-\cfrac{\omega_{0}}{z-\alpha_{1}-K_{\nu\,\boxright \,\mu}(z)}}
$$
which allows us to find an analytic form of $G_{\mu\,\boxplus\,\nu}(z)$ once we have 
an analytic formula for $K_{\nu\,\boxright \,\mu}(z)$. This type of algebraic computation 
was used, for instance, in [11], where we refer the reader for a general explicit formula for the 
Green function (equivalent to the Cauchy transform). Explicit computations based on this formula
for specific graphs are given in Sections 10-11.\\[10pt]
\myownsection
\begin{center}
{\sc 9. Quantum decomposition of adjacency matrices}
\end{center}
In this Section we will introduce a new type of `quantum decomposition' of the adjacency matrix $A({\cal G})$
of a graph ${\cal G}$ 
in which the distance is measured with respect to a set of vertices instead of a single vertex.

The results of this Section can be applied to any graph, not only free products of graphs to which 
this paper is devoted. In the latter case, it is of different category than those studied in 
the previous sections, but it bears some resemblance to the comb-orthogonal decomposition since
in some sense it `begins' with one copy of one of the graphs. However, as in the 
standard `quantum decomposition' [12,13], its components, called 
`quantum components', use infinitely many copies of both graphs and, moreover, 
cannot be represented as subgraphs of the product graph.
More importantly, it allows us to obtain more complete information 
on the spectral properties of many graphs, including certain free products, which 
cannot be obtained by means of other decompositions, including the standard `quantum decomposition'.

The set of vertices with respect to which distance is measured will be denoted by 
${\cal V}_{0}$. The sequence of sets
$$
{\cal V}_{n} = \{ v\in V; d(v,{\cal V}_{0})=n \},
$$
where $d(v,{\cal V}_{0}) = \min\{\,d(v,v_0)\,;\, v_0\in {\cal V}_{0}\,\}$
and $n\in {\mathbb N}^{*}:={\mathbb N}\cup \{0\}$ will be called the 
{\it distance partition} of the set $V$. 
The associated Hilbert space decomposition is of the form
\begin{equation}
l_{2}(V)=\bigoplus_{n\in {\mathbb N}^{*}}l_{2}({\cal V}_{n})
\end{equation}
which, in turn, leads to the {\it quantum decomposition} of $A$, by which
we shall understand the triple $(A^{+}, A^{0}, A^{-})$ of operators
on $l_{2}(V)$ given by
\begin{equation}
A^+ \delta(x)= \sum_{\stackrel{ y\sim x}{\scriptscriptstyle y\in {\cal V}_{n+1}}} \delta(y)\,,\ \ \ \
A^0 \delta(x)= \sum_{\stackrel{ y\sim x}{\scriptscriptstyle y\in {\cal V}_{n}}} \delta(y)\,,\ \ \ \
A^- \delta(x)= \sum_{\stackrel{ y\sim x}{\scriptscriptstyle y\in {\cal V}_{n-1}}} \delta(y)
\end{equation}
whenever $x\in {\cal V}_{n}$. Clearly, we have
$A=A^{+}+A^{0}+A^{-}$, which justifies the above terminology,
and, moreover, $(A^+)^*=A^-$ and $(A^0)^*=A^0$.
Finally, a non-zero vector $\xi\in l_2(V)$ will be called a 
{\it vacuum vector} if $A^-\xi = 0$. Of interest to us will be
vacuum vectors of special type.\\
\indent{\par}
{\sc Definition 9.1.}
A vector $\xi\in l_{2}(V)$ will be called a {\it J-vacuum vector} 
with respect to the quantum decomposition $(A^{+},A^{0}, A^{-})$ 
(or, simply, a J-vacuum vector) if it is a vacuum vector and 
for every $n\in {\mathbb N}^{*}$ it holds that
\begin{eqnarray}
A^-A^+(A^{+n}\xi)&=&\omega_n A^{+n}\xi\\
A^0(A^{+n}\xi)&=&\alpha_n A^{+n}\xi
\end{eqnarray}
where $\alpha_{n}\in {\mathbb R}$ and $\omega_{n}\geq 0$ and
where we use the convention that $A^{+m}\xi=0$ implies $\omega_m=\alpha_m=0$ and
thus $\omega_n=\alpha_n=0$ for $n\geq m$.\\
\indent{\par}
In this fashion we can associate with each J-vacuum vector
the Jacobi parameters written in the form of a pair of sequences $(\alpha, \omega)$ (called from now on {\it J-sequences}), where $\alpha =(\alpha_{n})_{n\in {\mathbb N}^{*}} $ and
$\omega =(\omega_{n})_{n\in {\mathbb N}^{*}}$.\\
\indent{\par}
{\sc Proposition 9.1.}
{\it J-vacuum vectors associated with non-identical J-sequences
are orthogonal.}\\
\indent{\par}
{\it Proof.}
Denote these vectors by $\xi$ and $\xi'$ and the associated J-sequences
by $(\alpha, \omega)$ and $(\alpha', \omega')$.
Assume that $\alpha_n \neq \alpha_n'$ for some $n\in \mathbb{N}^{*}$.
Without loss of generality we can assume that
$\alpha_n \neq 0$, which implies that
$\omega_{n-1}\neq 0$. Then
$$
\big< A^{+n}\xi,A^{+n}\xi' \big> =
\frac{1}{\alpha_n} \big< A^0 A^{+n}\xi,A^{+n}\xi' \big> =
\frac{\alpha_n'}{\alpha_n} \big<  A^{+n}\xi,A^{+n}\xi' \big>\,
$$
and therefore, $A^{+n}\xi \perp A^{+n}\xi'$, which gives
$$
0 = \big<  A^{+n}\xi,A^{+n}\xi' \big> = \big<  A^{-n}A^{+n}\xi,\xi' \big> = \omega \, \big<  \xi, \xi' \big>,
$$
where $\omega = \omega_{0}\omega_{1}\ldots \omega_{n-1}\neq 0$ and thus $\xi \perp \xi'$.
Similar computations for the case when $\omega_n \neq \omega_n'$
for given $n\in \mathbb{N}$ lead to orthogonality $\xi\perp \xi'$ as well.$\hfill{\blacksquare}$\\
\indent{\par}
For a given distance partition of $V$, any set $\Xi$ of vectors from $l_{2}(V)$ will be called {\it distance-adapted} if $\Xi=\bigcup_{n\in {\mathbb N}^{*}}\Xi_{n}$, where $\Xi_{n}\subseteq l_{2}({\cal V}_{n})$. Such sets are convenient to deal with and for that reason we 
show in the proposition given below that if we have a set of mutually orthogonal 
J-vacuum vectors, which we call an {\it orthogonal J-vacuum set}, 
we can always choose one which is distance-adapted.
For a given quantum decomposition of $A$, we denote by
$[\mathcal{A} x]$ and $[\mathcal{A}^+ x]$ the closed linear subspaces
generated by vectors $\{A^n x\}_{n=0}^{\infty}$ and
$\{A^{+n} x\}_{n=0}^{\infty}$, respectively.
\\
\indent{\par}
{\sc Proposition 9.2.}
{\it For a given distance partition of $V$ and the associated quantum decomposition of $A$,
let $\Xi$ be an orthogonal J-vacuum set. Then there exists an orthogonal
J-vacuum set $\Theta$, which is distance-adapted and such that}
$$
\bigoplus_{\xi\in\Xi} \; [\mathcal{A}\xi] \ \subseteq \
\bigoplus_{\xi\in\Theta}\ [\mathcal{A}\xi]
$$.
\indent{\par}
{\it Proof.}
Let $\Xi = \{ \,\xi_i\, ; i\in I\}$, for a countable set of indices $I$.
According to the decomposition (8.1), we have
$$
\xi_i = \sum_{n=0}^{\infty}\xi_i^{(n)}
$$
for every $i\in I$, where $\xi_i^{(n)} \in l_2({\cal V}_n)$ and $n\in {\mathbb N}^{*}$.
It is not difficult to observe that each $\xi_i^{(n)}$ is a J-vacuum  vector.
For every $n\in {\mathbb N}^{*}$, we choose from the set $\{\xi_i^{(n)} ; i\in I\}$ 
a maximal linearly independent set, which we denote $\Gamma_n = \{\gamma_1,\ldots,\gamma_{k_n}\}$.
Of course, $k_n \leqslant |{\cal V}_n|$. We divide $\Gamma_n$ into disjoint classes
$$
\Gamma_n = \Gamma_n(1) \cup \Gamma_n(2) \cup \ldots \cup \Gamma_n(l_n)
$$
subject to the condition: $\gamma_i,\gamma_j \in \Gamma_n(l) \iff \gamma_i,\gamma_j$
are associated with the same J-sequences. 
From Proposition 8.1 it follows that vectors from different classes are orthogonal.
Let  $\Theta_n(l)$ be the set obtained by
applying the Gram-Schmidt orthogonalization to class $\Gamma_n(l)$. Then
$$
\Theta_n = \Theta_n(1) \cup \Theta_n(2) \cup \ldots \cup \Theta_n(l_n)
$$
is an orthogonal set. Moreover, each element of $\Theta_n$
is a J-vacuum vector since it is a linear combination of J-vacuum vectors
associated with the same J-sequence. Finally, we take $\Theta$ to be the union of the $\Theta_n$, i.e.
$\Theta = \bigcup_{n\in {\mathbb N}^{*}} \Theta_n$.
From the construction of the set $\Theta$ it follows easily that
it satisfies the conditions stated above.$\hfill{\blacksquare}$\\
\indent{\par}
{\sc Definition 9.2.}
For a given distance-adapted J-vacuum set $\Xi$ with decomposition
$\Xi=\bigcup_{n\in {\mathbb N}^{*}}\Xi_{n}$,
define a sequence of mutually orthogonal sets by the recurrence
\begin{equation}
B_{0} = \Xi_0\,,\ \ \ \ \ B_{n+1} =
( A^+ B_n \cup \Xi_{n+1} ) \backslash \{0\}
\end{equation}
The set $\Xi$ will be called {\it generating} if
for every $n\in {\mathbb N}^{*}$ the set  $B_{n}$ is a basis in $l_2({\cal V}_n)$.
This notion will turn out useful in the theorem given below.
\\
\indent{\par}
{\sc Theorem 9.3.}
{\it If $\;\Xi\;$ is an orthogonal J-vacuum set
which is generating and distance-adapted, then we have the direct sum decomposition
$l_2(V)\; =\; \bigoplus_{\xi\in\Xi} \; [\mathcal{A}\xi]$}.\\
\indent{\par}
{\it Proof.}
First, let us show that for $\xi\in\Xi$ it holds that
$[\mathcal{A} \xi] = [\mathcal{A}^+ \xi]$. Notice that for $m < n$ we have
\begin{eqnarray*}
\left< A^{+n}\xi,A^{+m}\xi \right> &=& \big< A^{+(n-m-1)}\xi,A^{-(m+1)}A^{+m}\xi \big> \\
 &=& \omega \big< A^{+(n-m-1)}\xi, A^{-}\xi \big> \ = \ 0\,,
\end{eqnarray*}
and therefore $\{A^{+n} \xi\}_{n=0}^{\infty}$ is an orthogonal set.
Applying the Gram-Schmidt orthogonalization to the set
$\{A^n \xi\}_{n=0}^{\infty}$, we obtain $\{A^{+n} \xi\}_{n=0}^{\infty}$,
and thus $[\mathcal{A} \xi] = [\mathcal{A}^+ \xi]$.
Next, observe that for different $\xi,\xi'\in\Xi$ and arbitrary
$m,n\geqslant 0$, the vectors $A^{+n}\xi$, $A^{+m}\xi'$ are orthogonal,
which follows from a straightforward induction.
Therefore, $[\mathcal{A}^+ \xi] \perp [\mathcal{A}^+ \xi']$.
Finally, we know by assumption that $B_n$ is a basis in $l_2({\cal V}_n)$, and therefore
$$
l_2(V) \;=\; \bigoplus_{n=0}^{\infty} \; l_2(V_n) \;=\; \bigoplus_{n=0}^{\infty}
\text{span} (B_n) \;\subseteq\; \bigoplus_{\xi\in\Xi} \; [\mathcal{A^+}\xi] \;
= \; \bigoplus_{\xi\in\Xi} \; [\mathcal{A}\xi]\,,
$$
which completes the proof since the reverse implication is obvious.$\hfill{\blacksquare}$\\
\indent{\par}
Let us observe that Theorem 9.3 gives an interacting Fock space decomposition 
of $l_{2}(V)$ since $[\mathcal{A}\xi]$ is an interacting Fock space for each $\xi\in \Xi$, in which
the set $\{A^{+n} \xi\}_{n=0}^{\infty}$ is a basis. 
The results of this Section give us sufficient conditions for an orthogonal 
decomposition of $l_{2}(V)$ of Theorem 9.3 to exist and that in turn allows 
us to get detailed information about the spectral properties of the adjacency matrix $A$, including spectral distributions associated with all vacuum vectors $\xi\in \Xi$
which appear in this decomposition. In particular, this also gives the spectrum of the considered graph. 
\\[10pt]

\myownsection
\begin{center}
{\sc 10. Trees}
\end{center}
In this Section we use the theory of Section 9 to 
consider the simplest examples of $n$-ary trees and homogenous trees.
Here, the set ${\cal V}_{0}$ will consist of one root, which considerably simplifies
the spectral analysis and the corresponding quantum decomposition agrees with that
used in the approach of Hora and Obata [12,13]].
However, our analysis goes a little further since we study spectral 
distributions associated with all cyclic vectors. \\
\begin{center}
{\it $n$-ary trees ${\mathbb T}_n$}
\end{center}
By Theorem 9.3, it suffices to find a distance-adapted generating J-vacuum set. Denote by $W_{n}$ the set of words in $n$ letters $a_1,a_2, \ldots , a_n$,
including the empty word.
Then there exists a natural bijection between $W_{n}$ and $V({\mathbb T}_n)$.
If $w$ labels a vertex of ${\mathbb T}_n$ for which $d(w,e)=m$,
then $a_1w,a_2w, \ldots ,a_nw$ denote `sons' of $w$.
Let $\Xi_{0}=\{\delta(e)\}$ and
$$
\Xi_{m}=\{ \sum_{j=1}^{k}\left(\delta(a_jw)-\delta(a_{k+1}w)\right), \;\;
1\leq k \leq n-1, \;w\in W_{n}, \; |w|=m-1\}
$$
for every natural $m\geq 1$. Clearly,
$\Xi =\bigcup_{m=0}^{\infty} \Xi_{m}$
is a distance-adapted orthogonal J-vacuum set. In fact, it is easy to see that
it is an orthogonal set of vacuum vectors. 
Next, observe that the cardinalities of sets
$$
{\cal B}_{m}=A^{+}({\cal B}_{m-1})\cup \Xi_{m}, \;\; m\geq 1
$$
with $B_{0}=\Xi_{0}$ satisfy the recurrence
$$
|B_{m}|=|B_{m-1}|+(n-1)n^{m-1}
$$
since $|\Xi_{m}|=n^{m-1}$, which gives $|B_{m}|=n^{m}={\rm dim}{\cal V}_m$
and thus $B_{m}$ is a basis of ${\cal V}_m$. By Theorem 8.3, we have a direct
sum decomposition of $l_{2}(V)$ into the sum of $[{\cal A}\xi]$, $\xi\in \Xi$.
Finally,
every $\delta(w)$ (and thus every $\xi\in\Xi$) 
is an eigenvector of both $A^{-}A^{+}$ and $A^{0}$ with 
eigenvalues $n$ and $0$, respectively, for every $w\in V(T_{m})$. 
Therefore, every $\xi$ is a J-vacuum vector and the associated 
J-sequences are $\omega_{k}(\xi)=n$ and $\alpha_{k}(\xi)=0$ for every $k$.
This shows that to every $\xi$ corresponds the same Cauchy transform, namely that
of the form (6.3) and thus the measure (6.4). Thus, the spectrum of
${\mathbb T}_n$ agrees with the support of that measure, which is the interval $[-2\sqrt{n},2\sqrt{n}]$.\\
\begin{center}
{\it Homogenous trees ${\mathbb H}_n$}
\end{center}
A similar approach can be used for homogenous trees.
Since ${\mathbb H}_n$ is a symmetric graph, any vertex can be chosen to be the root denoted
$e$.  

If $n=2$, then $\Xi=\{\delta(e)\}$ and the decomposition of Theorem 9.3 
is simply $l_{2}(V)=[{\cal A}\delta(e)]$.  
Therefore, let $n\geq 3$ and denote by $a_1,a_2, \ldots , a_n$ the 'sons' of $e$ and
label all the vertices with distance bigger than $2$ from the root using
only $a_1,a_2, \ldots , a_{n-1}$. The situation is very similar to that
of the $n$-ary trees. Let
$\Xi_{0}=\{\delta(e)\}$ and
\begin{eqnarray*}
\Xi_{1}&=&\{\sum_{j=1}^{k}\left(\delta(a_j)-\delta(a_{k+1})\right), \; 1 \leq k \leq n-1\}\\
\Xi_{m}&=&\{\sum_{j=1}^{k}\left(\delta(a_jw)-\delta(a_{k+1}w)\right), \; 1 \leq k \leq n-2, w\in W_{n-1}\}
\end{eqnarray*}
for $m\geq 2$, where $W_{n-1}$ is the set of words in
$a_1, a_2, \ldots ,a_{n-1}$. 
The set $\Xi=\bigcup_{m=0}^{\infty}\Xi_{m}$ 
is clearly  a distance-adapted orthogonal set of vacuum vectors.
To show that $\Xi$ is generating, take the sequence $(B_{m})$ defined by
(7.5) and observe that we have the recurrence
$$
|B_{m}|=|B_{m-1}| + n(n-1)^{m-1}(n-2)
$$
since $|\Xi_{m}|=n(n-1)^{m-1}(n-2)$,
which is solved by $|B_{m}|=n(n-1)^{m-1}={\rm dim}{\cal V}_{m}$, and thus
$B=\bigcup_{m\in {\mathbb N}^{*}}B_{m}$ is a basis in $l_{2}(V)$.

The vector $\delta(e)$ (and thus $\xi_{0}$)
is an eigenvector of both $A^{-}A^{+}$ and $A^{0}$ with eigenvalues $n$ and 
$0$, respectively. This gives J-sequences
$\omega_0(\xi_0) = n$ and $\omega_{k}(\xi_0)=n-1$ for $k\geqslant 1$
with $\alpha_{k}(\xi)=0$ for all $k$, and thus the Cauchy transform 
$$
G_{\xi_0}(z) = \frac{-2z + nz - n{\sqrt{4 - 4n + z^2}}}{2\left( n - z \right)\left( n + z \right) }.
$$
If $n=1$ the measure $\mu_{\xi_0}$ consists of two atoms at $\pm 1$ of mass 1/2 each.
For $n\geqslant 2$, $\mu_{\xi_0}$ is absolutely continuous w.r.t the Lebesgue measure and
has density
$$
f_{\xi_0}(x) = \left| \frac{n\sqrt{4n-4 - x^2}}{2\pi\left( n - x \right)\left( n + x \right)} \right|
$$
on the interval $[-2\sqrt{n-1},2\sqrt{n-1}]$.

In turn, vectors $\delta(w)\neq \delta(e)$ (and thus $\xi\neq \xi_{0}$) are
eigenvectors of both $A^{-}A^{+}$ and $A^{0}$ with eigenvaluse $n-1$ and 
$0$, respectively.
Therefore, $\xi\neq \xi_{0}$ (for $n\geqslant 2$) are 
J-vacuum vectors with the Jacobi parameters of the form
$\omega_k(\xi) = n-1$ and $\alpha_{k}(\xi)=0$ for all $k$, which give the Cauchy transform
$$
G_{\xi}(z) = \frac{z- {\sqrt{z^2-4(n-1)}}}{2( n - 1)}.
$$
and thus $\mu_{\xi}$ is the Wigner measure on the interval $[-2\sqrt{n-1},2\sqrt{n-1}]$.
Although we get two different spectral distributions, their contributions to
the spectrum ${\rm spec}({\mathbb H}_n)=[-2\sqrt{n-1},2\sqrt{n-1}]$ are the same. \\[10pt]

\myownsection
\begin{center}
{\sc 11. Other examples}
\end{center}
In this Section we apply the same approach as in Section 10 to 
two types of examples: free products of complete graphs ${\mathbb K}_n*{\mathbb K}_m$ 
and free products of a complete graph with a graph of `fork'-type
${\mathbb K}_n*{\mathbb F}_m$. In the first case, we deal with free products of 
symmetric graphs and therefore their spectra can be determined 
from one spectral distribution (see [11] for a detailed study).
However, since our approach examines spectral distributions associated
with all cyclic vectors, even in this case we obtain new information 
(for instance, about the `multiplicity' of the spectrum).\\
\begin{center}
{\it Free products ${\mathbb K}_n*{\mathbb K}_m$} 
\end{center}
By $K_n$ we denote the {\it complete graph} with $n+1$ vertices, i.e.
such in which each pair of vertices forms an edge.
Observe that any choice of a root gives an isomorphic rooted graph since 
all vertices are equivalent.
Take two complete graphs ${\mathbb K}_n$ and ${\mathbb K}_m$ with vertices $x_0,x_1,\ldots,x_n$ (choose $x_{0}$ 
as the root) and $y_0,y_1,\ldots,y_m$ (choose $y_{0}$ as the root), 
respectively, where $n,m\geqslant 1$. 
Then ${\mathbb K}_n * {\mathbb K}_m = (V,E,e)$, where the distance partition of the
set of vertices $V$ can be given by the recursion
\begin{eqnarray}
{\cal V}_0 &=& \{e,x_1, \ldots ,x_n\} \nonumber \\
{\cal V}_{k+1} &=& \{ y_i w ;\ w\in {\cal V}_{2k},\ i=1,\ldots,m \}\\
{\cal V}_{2k+2} &=& \{ x_i w ;\ w\in {\cal V}_{2k+1},\ i=1,\ldots,n \} \,,\nonumber
\end{eqnarray}
where $k=0,1,\ldots$ (see Fig.1 and Fig.2). 
Directly from this construction it follows that 
\begin{equation}
|{\cal V}_{2k}|=(n+1) m^k n^k \ \ \ \text{and} \ \ \ |{\cal V}_{2k+1}|=(n+1) m^{k+1} n^k\,.
\end{equation}

\begin{picture}(145.00,60.00)(-10.00,-30.00)
\unitlength=1pt

\put(-25,-105){\includegraphics[]{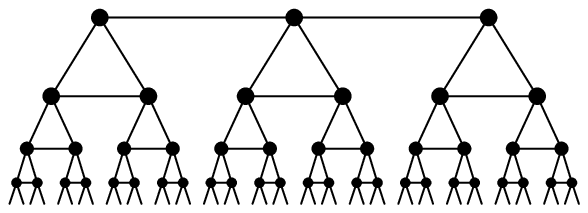}}

\put(190,-52){\includegraphics[]{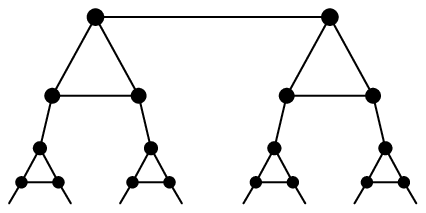}}

\bezier{500}(39,17)(95,35)(151,17)



\put(335,17){$\Big\} {\cal V}_0$}
\put(335,-7){$\big\} {\cal V}_1$}
\put(335.5,-22){$\} {\cal V}_2$}
\put(335.5,-36){$\} {\cal V}_3$}
\put(336,-53){$\vdots$}

\put(60,-65){{\it Fig. 6.} ${\mathbb K}_2 * {\mathbb K}_2$.}
\put(230,-65){{\it Fig. 7.} ${\mathbb K}_1 * {\mathbb K}_2$.}

\end{picture}

Let $(A^+,A^-,A^0)$ be the quantum decomposition of the incidence matrix of 
${\mathbb K}_n * {\mathbb K}_m$. We shall find the corresponding orthogonal 
J-vacuum set $\Xi$ which is generating and distance-adapted. 
Let
\begin{eqnarray}
\Xi_0 &=& \big\{\sum_{j=0}^n \delta(x_j)\big\} \cup \big\{\sum_{j=0}^{i-1} (\delta(x_j)-\delta(x_i));\ i=1,\ldots,n\big\} \nonumber \\
\Xi_{2k+1} &=& \big\{ \sum_{j=1}^{i-1} \left(\delta(y_j w)-\delta(y_i w)\right);\ w\in {\cal V}_{2k},\   i=2,\ldots,m \big\}\\
\Xi_{2k+2} &=& \big\{ \sum_{j=1}^{i-1} \left(\delta(x_j w)-\delta(x_i w)\right);\ w\in {\cal V}_{2k+1},\ i=2,\ldots,n \big\}\,,\nonumber
\end{eqnarray}
where $x_0$ is identified with $e$. 
This gives
$$
|\Xi_{2k+1}|=(m-1)|{\cal V}_{2k}|\:\;\;{\rm and}\;\;\; |\Xi_{2k+2}|=(n-1)|{\cal V}_{2k-1}|
$$
for $k\in {\mathbb N}^{*}$, which, by (9.2) leads to
\begin{equation}
|\Xi_{k}|=|{\cal V}_{k}|-|{\cal V}_{k-1}|.
\end{equation}
It is an elementary computation which 
shows that $\Xi_{k}$ is an orthogonal set for every $k\in {\mathbb N}^{*}$. 
Moreover, it is clear that these sets are mutually orthogonal and that they contain 
only vacuum vectors 
(use $A^{-}\delta(y_{j}w)=\delta(w)$ and $A^{-}\delta(x_{j}w')=\delta(w')$)
and therefore $\Xi$ is an orthogonal set of vacuum vectors.

We will show now that it is an orthogonal J-vacuum set.
Let $\xi\in\Xi_{2k+1}$ and $r\in\mathbb{N}$. We have
$$
A^{+r}\xi = A^{+r} \sum_{j=1}^{i-1} \left(\delta(y_jw)-\delta(y_{i}w)\right) = \sum_{j=1}^{i-1} \sum_{\stackrel{ |u|=r} {\scriptscriptstyle uy_{i}w \in V}} \left(\delta(u y_jw)-\delta(u y_{i}w)\right)
$$
for some $w\in {\cal V}_{2k}$. If $r$ is odd, then each $u$ in the above sum 
begins with a vertex from ${\mathbb K}_n$ and therefore each $\delta(uy_{k}w)$ 
in the above sum (and thus also $A^{+ r}\xi$) 
is an eigenvector of $A^{-}A^{+}$ with eigenvalue $m$. 
Moreover, each $\delta(uy_{k}w)$ (and thus also $A^{+r}\xi$) is an 
eigenvector of $A^{0}$ with eigenvalue $n-1$. 
Therefore, $\omega_{r}(\xi)=m$ and $\alpha_{r}(\xi)=n-1$ for $r$ odd.
In the case of $r$ even, each $u$ begins with a vertex from ${\mathbb K}_m$ in the above
reasoning and thus $\omega_{r}(\xi)=n$ and $\alpha_{r}(\xi)=m-1$. Analogous
computations holds when $\xi\in \Xi_{2k}$. Therefore, all vectors from $\Xi$
are J-vacuum vectors.

To show that $\Xi$ is generating, we need to check that for every $k\in{\mathbb N}^{*}$ the set $B_k$, defined recursively by  (9.5), is a basis in $l_2({\cal V}_{k})$. For that purpose it is enough to show that 
$|B_k|=|{\cal V}_k|$ since $B_k$ is a set of non-zero orthogonal vectors.
This is shown by induction.
For $k=0$ we have $|B_0| = |\Xi_0| = n+1 = |{\cal V}_{0}|$. Assume now that  $|B_k|=|{\cal V}_k|$ for some $k\in\mathbb{N}$. Then from (9.4) we get
$$
|B_{k+1}| = |A^+ B_{k}| + |\Xi_{k+1}| = |{\cal V}_k| + |{\cal V}_{k+1}| - |{\cal V}_{k}| = 
|{\cal V}_{k+1}|\,
$$
and our claim holds by induction. This proves that $\Xi$ is a generating 
J-vacuum set.

In order to find the spectrum of $A$, we find 
the probability measures $\mu_{\xi}$ with moments $\mu_{\xi}(n) = \left< A^n\xi, \xi \right>$
for every $\xi\in \Xi$ and then take the union of their supports. Let
$$
\xi_0 = \delta(e) + \delta(x_1) + \ldots + \delta(x_n)
$$
The associated J-sequences are of the form
$$
\omega_k(\xi_0) =
\begin{cases}
n &\text{$k$ odd}\\
m &\text{$k$ even}
\end{cases},\
\alpha_k(\xi_0) =
\begin{cases}
n   &\text{$k=0$}\\
m-1 &\text{$k$ odd}\\
n-1 &\text{$k$ even\; positive}\\
\end{cases}\\
$$
and they give a continued-fraction representation of the Cauchy transform $G_{\xi_0}$ of the measure $\mu_{\xi_0}$, which leads to
$$
G_{\xi_0}(z) = \frac{1 - mn + mz + nz - z^2 +
    {\sqrt{w(z)}}}{2\left( -1 + n - z \right)
    \left( m + n - z \right) }\,,
$$
$$
w(z) = \left( 1 - 2n + mn + \left( 2 - m - n \right) z + z^2 \right)^2-4m\left( 1 - m + z \right) \left( 1 - n + z \right).
$$
Applying the Stieltjes inversion formula to the transform $G_{\xi_0}$, we obtain the absolutely continuous part of $\mu_{\xi_0}$ of the form
$$
f_{\xi_0}(x) = \left| \frac{{\sqrt{-w(x)}}}{2\pi\left( x + 1 - n  \right)
    \left( x - m - n  \right)} \right|\,,\ \ \ \ \ \ \ \ x\in I_{m,n}\,,
$$
on the set $I_{m,n}$ being the union of two closed intervals
with ends at $\frac{1}{2} (m + n -2 \pm {\sqrt{4(\sqrt{m}\pm\sqrt{n})^2 + (m - n)^2}})$ 
and disjoint interiors. In addition, $\mu_{\xi_0}$ has an atom 
at $n-1$ of mass $\frac{1}{1+m}\max \{0,m-n\}$.

Let now $\xi\in\Xi_{2k}$ and $\xi\neq \xi_0$. Then the J-sequences
are of the form
$$
\omega_k(\xi) =
\begin{cases}
n &\text{$k$ odd}\\
m &\text{$k$ even}
\end{cases},\
\alpha_k(\xi) =
\begin{cases}
-1  &\text{$k=0$}\\
m-1 &\text{$k$ odd}\\
n-1 &\text{$k$ even \;positive}\\
\end{cases}\\
$$
which give the Cauchy transform 
$$
G_{\xi}(z) = \frac{1 - 2m + 2n - mn + (2-m+n)z+z^2+{\sqrt{w(z)}}}{2n\left( 1 - m + z \right) \left( 2 + z \right)}\,.
$$
Again, the Stieltjes inversion formula gives the absolutely continuous part
of $\mu_{\xi}$ of the form
$$
f_{\xi}(x) = \left| \frac{{\sqrt{-w(x)}}}{2\,n\,\left( 1 - m + x \right) \,\left( 2 + x \right) } \right|\,,\ \ \ \ \ \ \ x\in I_{m,n}\,.
$$
Besides, $\mu_{\xi}$ has atoms at $-2$ and $m-1$ of masses 
$\frac{mn-1}{n(1+m)}$ and $\frac{1}{n(1+m)}\max \{0,n-m\}$, respectively.

\vspace{2pt}
For $\xi\in\Xi_{2k+1}$, the J-sequences are of the form
$$
\omega_k(\xi) =
\begin{cases}
m &\text{$k$ odd}\\
n &\text{$k$ even}
\end{cases},\
\alpha_k(\xi) =
\begin{cases}
-1  &\text{$k=0$}\\
n-1 &\text{$k$ odd}\\
m-1 &\text{$k$ even\;positive}\\
\end{cases}\\
$$
which give the Cauchy transform 
$$
G_{\xi}(z) = \frac{1 - 2n + 2m - mn + (2-n+m)z+z^2+{\sqrt{w(z)}}}{2m\left( 1 - n + z \right) \left( 2 + z \right)}\,.
$$
$$
f_{\xi}(x) = \left| \frac{{\sqrt{-w(x)}}}{2m\left( 1 - n + x \right) \left( 2 + x \right) } \right|\,,\ \ \ \ \ \ \ x\in I_{m,n}\,.
$$
Besides, $\mu_{\xi}$ has atoms at $-2$ and $n-1$ of masses 
$\frac{mn-1}{m(1+n)}$ and $\frac{1}{m(1+n)}\max \{0,m-n\}$, respectively.

We conclude that for any $m,n\geqslant1$ the continuous part of the spectrum of  
${\mathbb K}_n * {\mathbb K}_m$ agrees with $I_{m,n}$. The point spectrum is given by
$$
\begin{cases}
\emptyset & \ \ \text{if} \ \ \ \ m=n=1\\
\{-2\}& \ \ \text{if} \ \ \ \ m=n>1\\
\{-2,m-1\}& \ \ \text{if} \ \ \ \ m<n\\
\{-2,n-1\}& \ \ \text{if} \ \ \ \ m>n
\end{cases}\,.
$$
\vspace{10pt}
\noindent
\begin{center}
{\it Free products ${\mathbb K}_n*{\mathbb F}_m$}
\end{center}
By ${\mathbb F}_m$ we denote the {\it fork of degree} $m$, i.e. a connected simple graph
with $m+1$ vertices, say $y_0,y_1, \ldots , y_m$, in which the only edges are
given by $y_{i}\sim y_{0}$ for $i\neq 0$.
As before, by $K_n$ we denote a complete graph with vertices $x_0,x_1,\ldots,x_n$. 
It is easy to see that the set of vertices of $K_n * {\mathbb F}_m $ coincides with 
that in the previous example. Besides, on $V$ we introduce 
the same distance partition as in (9.1) (see Fig.3). 
Let $\Xi = \sum_{i=0}^{\infty} \Xi_i$ be the set of vectors 
defined as in (11.3). That $\Xi$ is a generating J-vacuum set is shown
as in the case of ${\mathbb K}_n*{\mathbb K}_m$.

\begin{picture}(145.00,70.00)(-5.00,-40.00)
\unitlength=1pt

\put(55,-110){\includegraphics[]{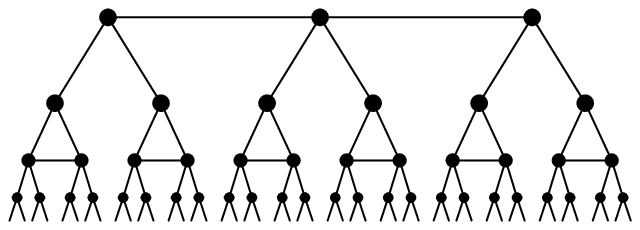}}

\bezier{500}(125,23)(186,43)(247,23)



\put(295,22){$\Big\} V_0$}
\put(295,-4){$\big\} V_1$}
\put(295.5,-20){$\} V_2$}
\put(295.5,-34){$\} V_3$}
\put(296,-53){$\vdots$}

\put(160,-70){{\it Fig 8.} ${\mathbb K}_2 * {\mathbb F}_2$.}

\end{picture}

Let us compute now the measures $\mu_{\xi}$ with moments $\mu_{\xi}(k)=\left< A^k\xi,\xi \right>$, for all $\xi\in\Xi$. The Jacobi parameters associated with vector $\xi_0$ are of the form
$$
\omega_k(\xi_0) = 
\begin{cases}
n &\text{$k$ odd}\\
m &\text{$k$ even}
\end{cases},\ 
\alpha_k(\xi_0) = 
\begin{cases}
n   &\text{$k=0$}\\
0   &\text{$k$ }\\
n-1 &\text{$k$ even}\\
\end{cases}\\
$$
The Cauchy transform $G_{\xi_0}$ of measure $\mu_{\xi_0}$ has the form
$$
G_{\xi_0}(z) = \frac{n - m - z - nz + z^2 - {\sqrt{ v(z) }}}{2 \left( m - n^2 + 2nz - z^2 \right) },
$$
where $v(z)={\left( m - n + z - nz + z^2 \right)}^2 - 4mz\left( z + 1 - n  \right)$. The measure $\mu_{\xi_0}$ is absolutely continuous w.r.t. the Lebesgue's measure and its density 
is of the form
$$
f_{\xi_0}(x) = \left| \frac{\sqrt{-v(x)}}{2\pi (m - n^2 + 2nx - x^2)} \right|, \ \ \ \ \ x\in J_{m,n}\,,
$$
where $J_{m,n}$ denotes the set, which is a union of two closed intervals with disjoint
interiors and ends at
$\frac{1}{2}(n-1 \pm \sqrt{4(\sqrt{m} \pm \sqrt{n})^2+(n-1)^2})$.

Let now $\xi\in\Xi_{2k}$ and $\xi\neq \xi_0$. Then the Jacobi parameters are of the form
$$
\omega_k(\xi) = 
\begin{cases}
n &\text{$k$ odd}\\
m &\text{$k$ even}
\end{cases},\ 
\alpha_k(\xi) = 
\begin{cases}
-1  &\text{$k=0$}\\
0   &\text{$k$ odd}\\
n-1 &\text{$k$ even\;positive}\\
\end{cases}\\
$$
which gives the Cauchy transform
$$
G_{\xi}(z) = \frac{ z^2 + z + n z -m + n - {\sqrt{ v(z) }}}{2 n \left( 1 - m + 2 z + z^2 \right) }.
$$
The absolutely continuous part of the measure $\mu_{\xi}$ is given by
$$
f_{\xi}(x) = \left| \frac{\sqrt{-v(x)}}{2\pi n (1 - m + 2 x + x^2)} \right|, \ \ \ \ \ x\in J_{m,n}\,.
$$
Besides, $\mu_{\xi}$ has two atoms at $\pm\sqrt{m}-1$ of mass  $\frac{n-1}{2n} $ each.

For $\xi\in\Xi_{2k+1}$, the Jacobi parameters are of the form
$$
\omega_k(\xi) = 
\begin{cases}
m &\text{$k$ odd}\\
n &\text{$k$ even}
\end{cases},\ 
\alpha_k(\xi) = 
\begin{cases}
n-1 &\text{$k$ odd}\\
0   &\text{$k$ even}\\
\end{cases}\\
$$
which gives the Cauchy transform
$$
G_{\xi}(z) = \frac{ m - n + z - nz + z^2 - {\sqrt{ v(z) }}}{2 m z}.
$$
The absolutely continuous part of $\mu_{\xi}$ is given by
$$
f_{\xi}(x) = \left| \frac{\sqrt{-v(x)}}{2\,\pi\, m\, z} \right|, \ \ \ \ \ x\in J_{m,n}\,.
$$
Besides, $\mu_{\xi}$ has an atom at $0$ of mass $\frac{1}{m}\max\{0,m-n\}$.

We conclude that the continuous spectrum of $K_n*{\mathbb F}_m$ coincides with $J_{m,n}$. 
Besides, $A$ has a discrete spectrum of the form
$$
\begin{cases}
\emptyset & \ \ \text{if} \ \ \ \ m=n=1\\
\{0\}& \ \ \text{if} \ \ \ \ m > n = 1\\
\{-2,0\}& \ \ \text{if} \ \ \ \ n > m = 1\\
\{\sqrt{m}-1,-\sqrt{m}-1\}& \ \ \text{if} \ \ \ \ n \geqslant m > 1\\
\{\sqrt{m}-1,-\sqrt{m}-1,0\}& \ \ \text{if} \ \ \ \ m > n > 1 
\end{cases}\,.
$$
\indent{\par}
For a general study of measures associated with mixed periodic Jacobi continued fractions, we
refer the reader to [14].\\

\begin{center}
{\sc Acknowledgements}
\end{center}
Two of the authors (R.L. and R.S.) would like to thank Professor Luigi Accardi for his hospitality during their
stay at the Volterra Center at the Universita di Roma Tor Vergata, which originated our joint work on this subject.\\

\newpage
\begin{center}
{\sc Bibliography}
\end{center}
[1] L. Accardi, A. Ben Ghorbal, N. Obata, Monotone independence, comb graphs and Bose-Einstein condensation,
{\it Infin. Dimens. Anal. Quantum Probab. Relat. Top.}, to appear.\\[3pt]
[2] K. Aomoto, Y. Kato, Green functions and spectra on free products of cyclic groups,
{\it Ann. Inst. Fourier} {\bf 38} (1988), 59-85.\\[3pt]
[3] D. Avitzour, Free products of $C^{*}$- algebras,
{\it Trans.~Amer.~Math.~Soc.} {\bf 271} (1982), 423-465.\\[3pt]
[4] Ph. Biane, Processes with free increments, {\it Math. Z.} {\bf 227} (1998), 143-174.\\[3pt]
[5] M. Bo\.{z}ejko, M. Leinert, R. Speicher, Convolution and limit theorems
for conditionally free random variables,
{\it Pacific J.Math.} {\bf 175} (1996), 357-388. \\[3pt]
[6] D.I. Cartwright, P.M. Soardi, Random walks on free products, quotients and amalgams,
{\it Nagoya J. Math.} {\bf 102} (1986), 163-180.\\[3pt]
[7] R. Burioni, D. Cassi, M. Rasetti, P. Sodano and A. Vezzani:
Bose-Einstein condensation on inhomogenuous complex networks,
{\it J. Phys. B: At. Mol. Opt. Phys.} \textbf{34} (2001), 4697--4710.\\[3pt]
[8] G. Giusiano, F.P. Mancini, P. Sodano, A. Trombettoni,
Local density of states in inhomogeneous graphs: the random walk and quantum probabilistic approaches,
Preprint, 16--11--2005.\\[3pt]
[9] A. Figa-Talamanca, T. Steger, Harmonic analysis for anisotropic
random walks on homogenous trees, {\it Memoirs Amer. Math. Soc.} {\bf 531} (1992).\\[3pt]
[10] C.D. Godsil, B. Mohar, Walk generating functions and spectral measures of infinite
graphs, {\it Linear Algebra Appl.} {\bf 107}(1988), 91-206.\\[3pt]
[11] E. Gutkin, Green's functions of free products of operators with applications to graph
spectra and to random walks, {\it Nagoya Math. J.} {\bf 149} (1998), 93-116.\\[3pt]
[12] A. Hora, N. Obata, {\it Quantum Probability and Spectral Analysis on Graphs},
monograph, to appear.\\[3pt]
[13] A. Hora, N. Obata, Quantum decomposition and quantum central limit theorem, {\it in}
"Fundamental Problems in Quantum Physics" (L. Accardi, S. Tasaki, Eds.),
pp. 284-305, World Scientfic, 2003.\\[3pt]
[14] Y. Kato, Mixed periodic Jacobi continued fractions, {\it Nagoya J. Math.} {\bf 104}
(1986), 129-148.\\[3pt]
[15] H. Kesten, Symmetric random walks on groups, {\it Trans. Am. Math. Soc.} {\bf 92}
(1959), 336-354.\\[3pt]
[16] G. Kuhn, Random walks on free products, {\it Ann. Inst. Fourier} {\bf 4} (1991), 467-491.\\[3pt]
[17] R. Lenczewski, private communication.\\[3pt]
[18] R. Lenczewski, Unification of independence in quantum probability,
{\it Infin. Dimens. Anal. Quantum Probab. Relat. Top.}, {\bf 1} (1998), 383-405.\\[3pt]
[19] R. Lenczewski, Reduction of free independence to tensor independence,
{\it Infin. Dimens. Anal. Quantum Probab. Relat. Top.} {\bf 7} (2004), 337-360.\\[3pt]
[20] R. Lenczewski, On noncommutative independence, in {\it QP-PQ: Quantum Probability
and White Noise Analysis. Vol. XVIII. Quantum Probability and Infinite Dimensional Analysis},
Eds. M. Schurmann, U. Franz, World Scientific, 2005, p. 320-336.\\[3pt]
[21] R. Lenczewski, Decompositions of the free additive convolution, preprint, archiv: math.OA/0608236,
2006.\\[3pt]
[22] Y.G. Lu, On the interacting Free Fock Space and the deformed Wigner Law,
{\it Nagoya Math. J.} 145 (1997) 1-28.\\[3pt]
[23] T. Matsui, BEC of Free Bosons on Networks,
{\it Infin. Dimens. Anal. Quantum Probab. Relat. Top.}, {\bf 1} (2006).\\[3pt]
[24] J.C. McLaughlin, Random walks and convolution operators on free products, dissertation,
New York University, 1986.\\[3pt]
[25] B. Mohar, W. Woess, A survey on spectra of infinite graphs,
{\it Bull. London Math. Soc.} {\bf 21}(1989), 209-234.\\[3pt]
[26] N. Muraki, Monotonic independence, monotonic central limit theorem and monotonic law of small numbers,
{\it Infin. Dimens. Anal. Quantum Probab. Relat. Top.} {\bf 4} (2001), 39-58.\\[3pt]
[27] N. Muraki, Monotonic convolution and monotone Levy-Hincin formula, preprint, 2000\\[3pt]
[28] N. Muraki, The five independences as quasi-universal products,
{\it Infinite Dim. Anal. Quantum Prob.} {\bf 5} (2002) 113-134. \\[3pt]
[29] N. Obata, Quantum probabilistic approach to spectral analysis of star graphs,
{\it Interdiscip. Inform. Sci.} {\bf 10} (2004), 41-52.\\[3pt]
[30] G. Quenell, Combinatorics of free product graphs, {\it Contemp. Math.} {\bf 206}
(1994), 257-281.\\[3pt]
[31] M. Sch\"urmann, Non-commutative probability on algebraic structures,
in {\it Probability Measures on Groups and Related Structures XI},
Proc., Oberwolfach 1994, ed. H. Heyer, World Scientific (1995) 332-356. \\[3pt]
[32] R. Speicher, R. Woroudi, Boolean convolution, in {\it Free Probability
Theory}, Ed. D. Voiculescu, pp. 267-279, Fields Inst. Commun. Vol.12, AMS, 1997.\\[3pt]
[33] D. Voiculescu, Symmetries of some reduced free product $\mathcal{C}^*$-algebras,
Operator Algebras and their Connections with Topology and Ergodic Theory,
{\it Lecture Notes in Math.} 1132, Springer, Berlin, 1985, 556-588.\\[3pt]
[34] D. Voiculescu, Addition of certain non-commuting random variables,
{\it J. Funct. Anal.} {\bf 66} (1986), 323-246.\\[3pt]
[35] D. Voiculescu, Noncommutative random variables and spectral problems in free product
$C^{*}$-algebras, {\it Rocky Mountain J. Math.} {\bf 20} (1990), 263--283.\\[3pt]
[36] D. Voiculescu, The analogues of entropy and of Fisher's information measure in free probability theory, I,
{\it Commun. Math. Phys.} {\bf 155} (1993), 71-92.\\[3pt]
[37] W. Woess, {\it Nearest neighbour random walks on free products of discrete groups},
{\it Boll. Unione Mat. Ital.} {\bf 5-B} (1986), 961-982.\\[3pt]
[38] W. Woess, {\it Random walks on Infinite Graphs and Groups},
Cambridge University Press, Cambridge, 2000.\\[3pt]
[39] D.W. Znojko, Free products of nets and free symmetrizers of graphs, {\it Mat.Sb. (N.S.)}
{\bf 98} (1975), 518-537.
\end{document}